\batchmode
\documentclass[11pt]{article}

\usepackage{epsfig}
\usepackage{graphicx}
\usepackage{color}
\usepackage{mathtools}

\newtheorem{theorem}{Theorem}
\newtheorem{lemma}{Lemma}
\newtheorem{corollary}{Corollary}

\newcommand{\be}{\begin{equation}}
\newcommand{\ee}{\end{equation}}
\newcommand{\bea}{\begin{eqnarray}}
\newcommand{\eea}{\end{eqnarray}}
\newcommand{\beas}{\begin{eqnarray*}}
\newcommand{\eeas}{\end{eqnarray*}}
\newcommand{\ba}{\begin{array}}
\newcommand{\ea}{\end{array}}

\DeclarePairedDelimiter{\floor}{\lfloor}{\rfloor}
\DeclarePairedDelimiter{\ceil}{\lceil}{\rceil}

\definecolor{armygreen}{rgb}{0.29, 0.33, 0.13}

\newcommand{\real}{\mbox{$\mathbb{R}$}}
\newcommand{\Natural}{\mbox{$\mathrm{I\!N}$}}
\newcommand{\Grad}{\ensuremath{\nabla}}
\newcommand{\eps}{\ensuremath{\epsilon}}

\newcommand{\bfq}{\ensuremath{\mathbf{q}}}

\newcommand{\bft}{\ensuremath{\mathbf{t}}}

\newcommand{\bfv}{\ensuremath{\mathbf{v}}}
\newcommand{\bfw}{\ensuremath{\mathbf{w}}}
\newcommand{\bfx}{\ensuremath{\mathbf{x}}}
\newcommand{\bfy}{\ensuremath{\mathbf{y}}}
\newcommand{\bfz}{\ensuremath{\mathbf{z}}}
\newcommand{\bfI}{\ensuremath{\mathbf{I}}}

\newcommand{\bfF}{\ensuremath{\mathbf{F}}}

\def\XXint#1#2#3{{\setbox0=\hbox{$#1{#2#3}{\int}$}
     \vcenter{\hbox{$#2#3$}}\kern-.5\wd0}}

\newcommand{\alpot}{\ensuremath{\frac{\alpha}{2}}}

\newcommand{\mcF}{\ensuremath{\mathcal{F}}}

\newcommand{\mcH}{\ensuremath{\mathcal{H}}}

\newcommand{\mcL}{\ensuremath{\mathcal{L}}}

\newcommand{\mcV}{\ensuremath{\mathcal{V}}}

\newcommand{\wtilde}[1]{\ensuremath{\widetilde{#1}}}

\def\qed{\hbox{\vrule width 6pt height 6pt depth 0pt}}

\usepackage{amsmath}
\usepackage{amssymb}

\title{A variable diffusivity fractional Laplacian} 
\author{	
	V.J.~Ervin\thanks{School of Mathematical and Statistical Sciences,
	  Clemson University, Clemson, South Carolina 29634-0975, USA.
	  email: {\tt vjervin@clemson.edu}. } }

\date{\today}

\begin{document}
\maketitle

\begin{abstract}
In this paper we analyze the existence, uniqueness and regularity of the solution to the generalized, 
variable diffusivity,
fractional Laplace equation on the unit disk in $\real^{2}$. For $\alpha$ the order of the differential operator,
our results show that for the symmetric, positive definite, diffusivity
matrix, $K(\bfx)$, satisfying $\lambda_{m} \, \bfv^{T} \, \bfv \ \le \, \bfv^{T} \, K(\bfx) \, \bfv \ \le \ \lambda_{M} \, \bfv^{T} \, \bfv$, 
for all $\bfv \in \real^{2}$, $\bfx \in \Omega$, with 
$\lambda_{M} \, < \, \frac{\sqrt{\alpha \, (2 + \alpha)}}{(2 - \alpha)} \, \lambda_{m}$, the 
problem has a unique solution. The regularity of the solution is given in an appropriately weighted Sobolev space in 
terms of the regularity of the right hand side function and $K(\bfx)$.
\end{abstract}



\setcounter{equation}{0}
\setcounter{figure}{0}
\setcounter{table}{0}
\setcounter{theorem}{0}
\setcounter{lemma}{0}
\setcounter{corollary}{0}
\setcounter{definition}{0}
%

%
\section{Introduction}
\label{sec_intro}

The diffusion operator is ubiquitous as a central modeling tool across a vast field of physical applications.
Characteristic of the usual diffusion process is the exponential spatial decay of the modeled quantity.
However, in a number of applications the spatial decay has been observed to be algebraic (``heavy tailed'')
leading to such processes being described as exhibiting anomalous diffusion. Nonlocal operators has been
an active area of research as tools to model anomalous diffusion. Examples of application
areas which have used such nonlocal operators include image processing \cite{bua101, gil081},
contaminant transport in ground water flow \cite{ben001}, electromagnetic fluids \cite{mcc811},
viscoelasticity \cite{mai971}, turbulent flow \cite{mai971, shl871}, chaotic dynamics \cite{zas931},
and biology \cite{bue141}.

The most commonly used approach to modeling anomalous diffusion processes is to use the fractional
Laplace operator. (Unlike the Laplace operator, the fractional Laplace operator is a nonlocal operator.)
Diffusion processes modeled using the fractional Laplace operator do exhibit algebraic spatial decay.

Unlike diffusion derived from random walk processes, many diffusion processes are driven by a local 
imbalance of a modeled quantity. For example, Fourier Law of heat conduction 
(\textit{heat flux is proportional
to the thermal gradient}) and Fick's Law of diffusion (\textit{molar flux is proportional to the concentration
gradient}) both posit this property. The proportionality constant is related to properties of the underlying medium, 
for example 
the thermal conductivity in the case of Fourier's Law.

For usual diffusion the impact of the
local imbalance is local. In the case of anomalous diffusion the impact of the local imbalance is nonlocal. For the
usual diffusion of the quantity $w$, the local flux in $\real^{d}$ is modeled as $K(\bfx) \Grad w(\bfx)$, and the spatial diffusion
operator as $- \Grad \cdot K(\bfx) \Grad w(\bfx)$, where $K(\bfx) \in \real^{d \times d}$ is a parameter related to the medium.
(In the heat equation 
$K(\bfx)$ represents the thermal conductivity of the medium.) To date the fractional Laplace operator has been restricted to
a constant material parameter, and formulated as an operator applied to an unknown quantity rather than an operator 
applied to the gradient of an unknown quantity.

In order to include a parameter related to the medium in an anomalous diffusion model we recently proposed the
following model \cite{zhe232}. For $\bfq(\bfx)$ denoting flux, $1 < \alpha < 2$
\[
    \begin{array}{ccc}
    \mbox{usual diffusion}    &  \hspace{0.75in}  &   \mbox{anomalous diffusion}   \\
   \Grad \cdot \bfq(\bfx) \ = \ - \Grad \cdot K(\bfx) \Grad w(\bfx)  &      &   
     \Grad \cdot \bfq(\bfx) \ = \  - \Grad \cdot \left(- \Delta \right)^{\frac{\alpha - 2}{2}} K(\bfx) \Grad w(\bfx) \, ,    
    \end{array}
\]
where $\left(- \Delta \right)^{\frac{\alpha - 2}{2}}(\cdot)$ defined in \eqref{defRPo}, represents the Riesz potential operator.
(For suitably nice functions (see \cite[Proposition 1]{dyd171}) the fractional Laplace operator is the inverse of the Riesz potential 
operator.)

For $\bfI$ the identity matrix in $\real^{d \times d}$, $k \in \real$, and $f$ a sufficiently nice function, a Fourier transform
argument shows that
\[
    k \left(- \Delta \right)^{\alpot} f(\bfx) \ = \ - \Grad \cdot \left(- \Delta \right)^{\frac{\alpha - 2}{2}} k \bfI \,  \Grad f(\bfx) \, , \ \ 
    \bfx \in \real^{d} \, .
\]

For $\Omega$ denoting the unit disk in $\real^{2}$, in \cite{zhe232} we showed, that for $K(\bfx)$ a symmetric, positive
definite, constant matrix and $f \in H_{\alpot}^{s}(\Omega)$, the existence and uniqueness of 
$\wtilde{u} \in \omega^{\alpot} \otimes H_{\alpot}^{s + \alpha}(\Omega)$ (see Section \ref{ssec_H01} for the definition of the spaces)
satisfying
\be
\mcL(\wtilde{u})(\bfx) \, := \, - \Grad \cdot \left(- \Delta \right)^{\frac{\alpha - 2}{2}} K(\bfx) \Grad \wtilde{u}(\bfx) \ = \ 
f(\bfx) \, , \ \bfx \in \Omega \, , \ \ \mbox{ with } \ \  \wtilde{u}(\bfx) \, = \, 0 \, , \  \bfx \in \real^{2} \backslash \Omega \, .
\label{defLL}
\ee

A complete analysis of the generalized fractional Laplacian (with advection and reaction terms) on a bounded
interval in $\real^{1}$ was recently obtained in \cite{zhe231}. (In \cite{zhe231} the generalized fractional Laplacian is 
a special case ($r = \frac{1}{2}$) of the more general operator studied therein.)


For $\wtilde{u} \, = \, \omega^{\alpot} \, u(\bfx)$, $u \in H^{1}_{\alpot}(\Omega)$ and 
$\wtilde{v} \, = \, \omega^{\alpot} \, v(\bfx)$, $u \in H^{\alpha - 1}_{\alpot}(\Omega)$ the weak formulation of \eqref{defLL} studied in 
Section \ref{sec_exu} is 
\begin{align*}
B(\wtilde{u} , \wtilde{v}) & = \  \left( K(\bfx) \Grad \wtilde{u} \, , \, (-\Delta)^{\frac{\alpha - 2}{2}} \Grad \wtilde{v} \right) \ = \ 
 F(\wtilde{v}) \, ,  \\
 \Longleftrightarrow \ & \ \left( \omega^{\alpot - 1} \, K(\bfx) \, U \, , \, V \right) \ = \ 
 F(\wtilde{v}) \, , 
\end{align*}
where $U$, $V \in \big( L_{\alpot -1}^{2}(\Omega) \big)^{2}$. The difficulty in establishing the existence and uniqueness of the
solution $\wtilde{u}$ lies in the fact that the subspaces in $\big( L_{\alpot -1}^{2}(\Omega) \big)^{2}$ containing $K(\bfx) \, U$
and $V$ are not equal. This means that establishing the appropriate inf-sup condition is satisfied is challenging. Essentially one needs to 
show that the norm of the piece of $K(\bfx) \, U$ that does not lie in the subspace containing $V$ can be bounded by
the norm of the piece of $K(\bfx) \, U$ that does lie in the subspace containing $V$. This leads to a restriction on the
symmetric, positive definite, matrix $K(\bfx)$ involving the parameter $\alpha$. As the parameter $\alpha \rightarrow 2$ (and
the operator $\mcL(\cdot)$ becomes the usual variable coefficient diffusion operator) the restriction on $K(\bfx)$ vanishes.
(See Theorem \ref{thmexun}.)

This paper is organized as follows. 
Section \ref{sec_prelim} contains definitions, notation and some results used later in the paper. 
Within the function space framework introduced in Section  \ref{ssec_H01}, Section \ref{sec_map} presents mapping properties of
the operators involved. The well posedness of the solution $\wtilde{u}$ to
\eqref{defLL}, along with its regularity as a function of the regularity of $f$,  is established in Section \ref{sec_exu}, subject to
suitable assumptions on $K(\bfx)$. Some concluding remarks are then given in Section \ref{sec_conc}.

 \setcounter{equation}{0}
\setcounter{figure}{0}
\setcounter{table}{0}
\setcounter{theorem}{0}
\setcounter{lemma}{0}
\setcounter{corollary}{0}
\setcounter{definition}{0}

\section{Preliminaries}
 \label{sec_prelim}

In this section we present definitions, notation and some results used later in the paper.

\textbf{The Fractional Laplacian} \\
In $\real^{d}$, the (integral) fractional Laplacian of a function $u(\bfx)$, 
$\left( - \Delta \right)^{\alpot} u(\bfx)$, is defined as
\[
 \left( - \Delta \right)^{\alpot} u(\bfx) \ := \ \frac{1}{| \gamma_{d}(-\alpha)|} \, \lim_{\eps \rightarrow 0} 
 \int_{\real^{d} \backslash B(\bfx, \eps)} \, \frac{u(\bfx) \, - \, u(\bfy)}{| \bfx \, - \, \bfy |^{d + \alpha}} \, d\bfy \, , 
\]
%
where $\gamma_{d}(\alpha) \, := \, 2^{\alpha} \, \pi^{d/2} \, \Gamma(\alpot) / \Gamma( (d - \alpha)/2 )$, and 
$B(\bfx, \eps)$ denotes the ball centered at $\bfx$ with radius~$\eps$.

\textbf{Riesz Potential Operator} \\
Related to the integral fractional Laplacian is the Riesz potential operator, which is defined for $\alpha \in (0, d)$ by \cite{dyd171}
\be
   \left( - \Delta \right)^{- \alpot} u(\bfx) \ := \ \frac{1}{\gamma_{d}(\alpha)} \, \int_{\real^{d}} 
   \frac{u(\bfx - \bfy)}{| \bfy |^{d - \alpha}} \, d\bfy \, .
\label{defRPo}
\ee

For suitably nice functions (see \cite[Proposition 1]{dyd171}) the integral fractional Laplacian is the inverse of the Riesz potential 
operator.

The following lemma is used in obtaining the weak formulation of \eqref{defLL} studied in Section \ref{sec_exu}.
\begin{lemma} \label{mvRPO}
For $\Omega$ a bounded domain in $\real^{d}$ and functions $f(\bfx)$ and $g(\bfx)$ satisfying $supp(f) , \, supp(g) \, \subset \Omega$
we have
\[
 \int_{\Omega} \left( - \Delta \right)^{- \alpot} f(\bfx) \ g(\bfx) \, d\bfx \ = \ 
  \int_{\Omega} f(\bfx) \ \left( - \Delta \right)^{- \alpot} g(\bfx) \, d\bfx \, .
\]
\end{lemma}
\textbf{Proof}:
Using \eqref{defRPo},
\begin{align*}
& \int_{\Omega} \left( - \Delta \right)^{- \alpot} f(\bfx) \ g(\bfx) \, d\bfx
\ = \  \int_{\Omega} \ \int_{\real^{d}} \, \frac{f(\bfx \, - \, \bfy)}{| \bfy |^{d - \alpha}} d\bfy \  g(\bfx) \, d\bfx  \\
& \quad =  \  \int_{\Omega} \ \int_{\real^{d}} \, \frac{f(\bfz) \, g(\bfx)}{| \bfx \, - \, \bfz |^{d - \alpha}} d\bfz \  \, d\bfx 
 \ \ \mbox{(letting $\bfz \, = \, \bfx \, - \, \bfy$)} \\
& \quad =  \  \int_{\real^{d}} \ \int_{\real^{d}} \, \frac{f(\bfz) \, g(\bfx)}{| \bfx \, - \, \bfz |^{d - \alpha}} d\bfz \  \, d\bfx 
 \ \ \mbox{(as $supp(g) \subset \Omega$)} \\
& \quad =  \  \int_{\bfz \in \real^{d}} \ \int_{\bfx \in \real^{d}} \, \frac{f(\bfz) \, g(\bfx)}{| \bfx \, - \, \bfz |^{d - \alpha}} d\bfx \  \, d\bfz 
 \ \ \mbox{(reversing the order of integration)} \\
& \quad =  \  \int_{\bfz \in \Omega} \ \int_{\bft \in \real^{d}} \, \frac{f(\bfz) \, g(\bft \, - \, \bfz)}{| \bft |^{d - \alpha}} d\bft \  \, d\bfz
 \ \ \mbox{(using $supp(f) \subset \Omega$, and letting $\bft \ = \ \bfx \, - \, \bfz$)} \\
& \quad = \   \int_{\Omega} f(\bfx) \ \left( - \Delta \right)^{- \alpot} g(\bfx) \, d\bfx \, .
\end{align*}
\mbox{ } \hfill \qed

\subsection{Jacobi Polynomials}
\label{ssec_jpoly}
The Jacobi polynomials are defined in terms of the regularized Hypergeometric function $_{2}\bfF_{1}(\cdot)$ 
(see \cite{dyd171}) as 
\begin{align*}
P_{n}^{(a , b)}(t) &:= \ \frac{\Gamma(a \, + \, 1 \, + \, n)}{n !} \,
 _{2}\bfF_{1}\left( \begin{array}{c}
                          -n , \ 1 + a + b + n \\  a + 1  \end{array} \, \Big{\vert} \, \frac{1 - t}{2} \right)   \\ 
&= \  \frac{(-1)^{n} \, \Gamma(b \, + \, 1 \, + \, n)}{n !} \,
 _{2}\bfF_{1}\left( \begin{array}{c}
                          -n , \ 1 + a + b + n \\  b + 1  \end{array} \, \Big{\vert} \, \frac{1 + t}{2} \right)   \\ 
&= \ (-1)^{n} \frac{\Gamma(n + 1 + b)}{n ! \, \Gamma(n + 1 + a + b)} \, \sum_{j = 0}^{n} (-1)^{j} \, 2^{-j} \, 
\left( \begin{array}{c}
           n  \\ j  \end{array} \right) \frac{\Gamma(n + j + 1 + a + b)}{\Gamma(j + 1 + b)} \, (1 + t)^{j} \, .   
\end{align*}


In case $a, \, b > -1$ the Jacobi polynomials satisfy the following orthogonality property.
\begin{align*}
 & \int_{-1}^{1} (1 - t)^{a} (1 + t)^{b} \, P_{j}^{(a , b)}(t) \, P_{k}^{(a , b)}(t)  \, dt 
 \ = \
   \left\{ \begin{array}{ll} 
   0 , & k \ne j  \\
   |\| P_{j}^{(a , b)} |\|^{2}
   \, , & k = j  
    \end{array} \right.  \, ,  \\  
& \quad \quad \mbox{where } \  \ |\| P_{j}^{(a , b)} |\| \ = \
 \left( \frac{2^{(a + b + 1)}}{(2j \, + \, a \, + \, b \, + 1)} 
   \frac{\Gamma(j + a + 1) \, \Gamma(j + b + 1)}{\Gamma(j + 1) \, \Gamma(j + a + b + 1)}
   \right)^{1/2} \, .
\end{align*}

\subsection{Solid Harmonic Polynomials}
\label{ssec_shpoly}
The solid harmonic polynomials in $\real^{d}$ are the polynomials in $d$ variables which satisfy 
Laplace's equation.

In $\real^{2}$ the solid harmonic polynomials of degree $l$ can be conveniently written in polar
coordinates
as $\{ r^{l} \cos (l \varphi ) \ , \   r^{l} \sin (l \varphi ) \} $.

\subsection{The Function Space $L^{2}_{\beta}(\Omega)$}
\label{ssec_funspc}     
With respect to the Cartesian coordinate system, we denote a point $\bfx \in \real^{2}$ as $\bfx = (x , y)$, and
with respect to the polar coordinate system, $\bfx = (r , \varphi)$.

We let $\Omega$ denote the unit disk in $\real^{2}$, i.e., $\Omega \, := \, \{ \bfx = (x , y) \, : \, x^{2} + y^{2} < 1\} . $

For the weight function $w(\cdot)$, $w(\bfx) > 0 \, ,
\bfx \in \Omega$, the associated weighted $L^{2}$ function space is defined as
\[
 L^{2}_{w}(\Omega) \, := \, \{ f \, : \, \| f \|_{L^{2}_{w}} < \infty \} \, , \ \ \mbox{ where }
  \| f \|_{L^{2}_{w}}^{2} \, := \, \int_{\Omega} \, w(\bfx) \, ( f(\bfx) )^{2} \, d\Omega \, .
\]
Associated with $ L^{2}_{w}(\Omega) $ we have the inner product, defined for $f, \, g \, : \, \Omega \rightarrow \real^{n}$, by
\[
  (f , g)_{L^{2}_{w}} \, := \, \int_{\Omega} w(\bfx) \, f(\bfx) \cdot g(\bfx) \, d\Omega \, .
\]

The usual $L^{2}(\Omega)$ norm and inner product (corresponding to $w(\bfx) = 1$) are respectively denoted by
$\| f \|$, and $(f , g)$.

For $\bfx = (r , \varphi) \in \Omega$ , the weight functions
$ \omega^{\gamma} \, := \, (1 - r^{2})^{\gamma} $, for $\gamma \in \real$, play a central role in the analysis.
To simplify notation, for lower case Greek letters, $\alpha, \beta, \gamma \in \real$ we let
\[
 (f , g)_{\beta} \, := (f , g)_{L^{2}_{\omega^{\beta}}} \, = \, \int_{\Omega} \omega^{\beta}  \, f(\bfx) \cdot g(\bfx) \, d\Omega
 \ = \ \int_{\Omega}(1 - r^{2})^{\beta}  \, f(\bfx) \cdot g(\bfx) \, d\Omega \, ,
\]
%
\[  
 \| f \|_{\beta} \, := \, (f , f)_{\beta}^{1/2} \ \ \mbox{and } \ \ L^{2}_{\beta}(\Omega) \, := \, L^{2}_{\omega^{\beta}}(\Omega) \, .
\]

Also, for notation convenience, let $\rho := \, 2r^{2} - 1$. \\

A basis for 
$L^{2}_{\beta}(\Omega)$ is given by a product of the solid harmonic polynomials and Jacobi polynomials.

\subsubsection{Basis for $L^{2}_{\beta}(\Omega)$}
\label{sssec_R2}     
\[
\mbox{Let } \ \mcV_{0 , 1}(\bfx) := \ \frac{1}{2} \, , \ 
\mcV_{l , 1}(\bfx) \ := \ r^{l} \, \cos(l \varphi) \, , \  \ l = 1, 2, \ldots  \ \ \mbox{ and } \
\mcV_{l , -1}(\bfx) \ := \ r^{l} \, \sin(l \varphi) \, , \  \ l = 1, 2, \ldots  .
\]

We also use the following notation
\[
    \mcV_{l , \mu^{*}}(\bfx) \ = \ \left\{ \begin{array}{rl}
    \mcV_{l , -1}(\bfx) & \mbox{  if } \mu = 1 \, ,   \\
    \mcV_{l , 1}(\bfx) & \mbox{  if } \mu = -1 \, .  \end{array} \right.
\]

Additionally, for a linear operator $\mcF \cdot$,  we use $\mcF(\mcV_{l , \mu}(\bfx)) \ = \ (\pm) \mcV_{l , \sigma}(\bfx)$ to denote
\[
\mcF(\mcV_{l , \mu}(\bfx)) = \ (\pm) \mcV_{l , \sigma}(\bfx) \ =  \ \left\{ \begin{array}{rl}
   + \mcV_{l , \sigma}(\bfx) & \mbox{  if } \mu = 1 \, , \\
   - \mcV_{l , \sigma}(\bfx) & \mbox{  if } \mu = -1 \, .  \end{array} \right.
\]
For example,
\[
 \frac{\partial}{\partial \varphi} \mcV_{l , \mu}(\bfx) \ = \ (\mp) \mcV_{l , \mu^{*}}(\bfx)
\]

An orthogonal basis for $L^{2}_{\beta}(\Omega)$ is  \cite{li141, wun051} (recall $\rho \, = \, 2 r^{2} \, - \, 1$)
\be
 \left\{  \cup_{l = 0}^{\infty} \cup_{n = 0}^{\infty} \left\{ \mcV_{l , 1}(\bfx) \, P_{n}^{(\beta  ,  l)}(\rho) \right\} \right\}   \ \cup
  \left\{ \cup_{l = 1}^{\infty} \cup_{n = 0}^{\infty}  \left\{ \mcV_{l , -1}(\bfx) \, P_{n}^{(\beta  ,  l)}(\rho) \right\} \right\}    \, .
\label{bsR2}
\ee

For notation brevity  we denote the basis in \eqref{bsR2} as
\[
  \cup_{l = 0}^{\infty} \cup_{n = 0}^{\infty}  \cup_{\mu = {1 , -1}} \left\{ \mcV_{l , \mu}(\bfx) \, P_{n}^{(\beta  ,  l)}(\rho) \right\} 
\]
where we implicit assume that the terms   $\mcV_{0, -1}(\bfx) \, P_{n}^{(\beta  ,  l)}(\rho), \, n = 0, 1, \ldots$ are omitted 
from the set. 

We trivially extend the definitions of $\mcV_{l, \mu}(\bfx)$ and $P_{n}^{(\alpha  ,  \beta)}(t)$ to negative integer values for $l$ and $n$
by $\mcV_{l, \mu}(\bfx) = 0$ for $l \le -1$, and $P_{n}^{(\alpha  ,  \beta)}(t) = 0$ for $n \le -1$.

In \cite{zhe232} we used the notation $V_{l, \mu}(\bfx)$ to represent the harmonic polynomials, with $\mcV_{l, \mu}(\bfx) = V_{l, \mu}(\bfx)$,
except for $\mcV_{0, 1}(\bfx) = \, \frac{1}{2} V_{0, 1}(\bfx) = \frac{1}{2}$. This change of notation, together with the extension described in
the previous paragraph allows a more compact representation for the derivative formulas in
Theorems \ref{gradinR2}, \ref{figradinR2}, and Corollary \ref{gradwR2}.

\subsection{The function space $H_{\gamma}^{s}(\Omega)$}
\label{ssec_H01}
Weighted Sobolev spaces are particularly useful in the analysis of problems whose solutions have a known 
singular behavior. In these cases the use of weighted Sobolev spaces afford a much more succinct analysis
than using the (unweighted) usual Sobolev spaces. The fact that the solutions to the fractional Laplace equation
have a known singular behavior at the boundary of the domain \cite{aco171, ros141} motivates the use of the
following spaces.

Similar to Bab\v{u}ska and Guo \cite{bab011} and Guo and Wang \cite{guo041}, define the 
weighted Sobolev space $H_{\gamma}^{s}(\Omega)$, for $s \in \Natural$, as
\be
H_{\gamma}^{s}(\Omega) \, := \, \Big\{ f \in L^{2}_{\gamma}(\Omega) \, : \, 
| f |_{H_{\gamma}^{s}} \, := \, \left( \sum_{j = 0}^{s} \left( \begin{array}{c}   s  \\ j  \end{array} \right) \, 
\left\| \frac{\partial^{s} f(\bfx)}{\partial y^{j} \, \partial x^{s - j}} \right\|^{2}_{\gamma + s} \right)^{1/2} \, < \, \infty \Big\} \, .
\label{defHsp}
\ee
We associate with $H_{\gamma}^{s}(\Omega)$ the norm
\[
  \| f \|_{H_{\gamma}^{s}} \, := \, \left( \| f \|_{\gamma}^{2} \ + \  | f |_{H_{\gamma}^{s}}^{2} \right)^{1/2} \, .
\]

For $s > 0$, $s \not\in \Natural$,  $H_{\gamma}^{s}(\Omega)$ is defined using the K-method of interpolation. For
$s < 0$ the spaces are defined by (weighted) $L^{2}_{\gamma}$ duality.

Following \cite{erv191}, introduce the space $W^{n, \infty}_{w}(\Omega)$ and its associated norm, defined for 
$n \in \Natural_{0}$ as
\begin{align*}
W^{n, \infty}_{w}(\Omega) &:= \ \Big\{ f \, : \, (1 - r^{2})^{m/2} \frac{\partial^{m} f(\bfx)}{\partial y^{j} \, \partial  x^{m-j}} 
\in L^{\infty}(\Omega) \, , \ 
\mbox{for } j = 0, 1, \ldots, m, \  m = 0, 1, \ldots, n \Big\} \, ,   \\   
\| f \|_{W^{n, \infty}_{w}} &:= \ \max_{0 \le j \le m \le n} 
\left\{ \left\| (1 - r^{2})^{m/2} \frac{\partial^{m} f(\bfx)}{\partial y^{j} \, \partial  x^{m-j}} \right\|_{L^{\infty}} \right\} \, .  
\end{align*}

The subscript $w$ denotes the fact that $W^{n, \infty}_{w}(\Omega)$ is a weaker space than $W^{n, \infty}(\Omega)$
in that the derivatives of functions in $W^{n, \infty}_{w}(\Omega)$ may be unbounded at the boundary of $\Omega$.

For $f(\cdot) \in H_{\gamma}^{s}(\Omega)$ we have the following result for the product $k(\cdot) \, f(\cdot)$.
\begin{lemma} \label{prodlma}
For $\gamma > -1$, $0 \le s \le n \in \Natural_{0}$, and $k(\cdot) \in W^{n, \infty}_{w}(\Omega)$. Then, for 
$ f(\cdot) \in H^{s}_{\gamma}(\Omega)$
\be
  k(\cdot) \, f(\cdot) \in H^{s}_{\gamma}(\Omega) \   
  \mbox{ with }  \  \| k \, f \|_{H^{s}_{\gamma}(\Omega)} \, \le \, 
  \| k \|_{W^{n, \infty}_{w}(\Omega)} \, \| f \|_{H^{s}_{\gamma}(\Omega)} \,  .
\label{wsp3} 
\ee  
\end{lemma}
\textbf{Proof}: The proof follows in a similar manner to that of Lemma 5.2 in \cite{erv191}.  \\
\mbox{ } \hfill \qed

In \cite{dyd171}, Dyda, Kuznetsov, and Kwa\'{s}nicki showed that the basis functions for $L^{2}_{\alpot}(\Omega)$
are pseudo eigenfunctions for the fractional Laplace operator of order $\alpha$ on the unit disk (see \eqref{bceq0v}).
In view of this property, it is natural to seek a solution expressed as an infinite series of these basis functions. The 
appropriate space to study such functions is defined by the decay rate of their coefficients.

For $s \in \real$, the space $\wtilde{H}_{\gamma}^{s}(\Omega)$ is defined as
\be
\wtilde{H}_{\gamma}^{s}(\Omega) \ := \ \Big\{ u(\cdot) \in L^{2}_{\gamma}(\Omega) \, : \, 
\sum_{l, n, \mu} (n + 1)^{s} (n + l + 1)^{s} \, a_{l, n, \mu}^{2} \, 
\| \mcV_{l, \mu}(\bfx) P_{n}^{(\gamma , l)}(\rho) \|_{L^{2}_{\gamma}}^{2}  < \infty \Big\} \, ,
\label{defXs}
\ee
where
\[
u(\bfx) \ = \ \sum_{l, n, \mu} a_{l, n, \mu} \, \mcV_{l, \mu}(\bfx) P_{n}^{(\gamma , l)}(\rho) \, , 
\]
with 
\[
\| u \|_{\wtilde{H}_{\gamma}^{s}}  := \ \Big( \sum_{l, n, \mu} (n + 1)^{s} (n + l + 1)^{s} \, a_{l, n, \mu}^{2} \, 
\| \mcV_{l, \mu}(\bfx) P_{n}^{(\gamma , l)}(\rho) \|_{L^{2}_{\gamma}}^{2} \Big)^{1/2}.
\]

In \cite{erv241} we showed that the spaces $H_{\gamma}^{s}(\Omega)$ and $\wtilde{H}_{\gamma}^{s}(\Omega)$
are equivalent, and that their dual spaces, with respect to the $L^{2}_{\alpot}(\Omega)$ pivot space, are characterized by
$\wtilde{H}_{\gamma}^{-s}(\Omega)$.

Consequently, we use $H_{\gamma}^{s}(\Omega)$ to denote the spaces defined in \eqref{defHsp} and \eqref{defXs}.

\subsection{Some useful derivative properties}
\label{ssec_derP}
In this section we present some derivative formulas used in the subsequent analysis.

\begin{theorem} \cite{zhe232} \label{gradinR2}
 Let $f(\bfx) \ = \ (1 - r^{2})^{\frac{\alpha}{2}} \mcV_{l , \mu}(\bfx) P_{n}^{(\frac{\alpha}{2} , l)}(2r^{2} - 1)$. Then,
 \begin{align}
\frac{\partial f}{\partial x} &= \
 - (n + \alpot) (1 - r^{2})^{\alpot - 1} \mcV_{l+1 , \mu}(\bfx) P_{n}^{(\alpot-1 , l+1)}(2r^{2} - 1) \nonumber \\
 & \mbox{ } \quad \quad \ - \
(n + 1) (1 - r^{2})^{\alpot - 1} \mcV_{l-1 , \mu}(\bfx) P_{n+1}^{(\alpot-1 , l-1)}(2r^{2} - 1) \, ,    \label{ders1} \\
\frac{\partial f}{\partial y} &= \
 - ( \pm ) (n + \alpot) (1 - r^{2})^{\alpot - 1} \mcV_{l+1 , \mu^{*}}(\bfx) P_{n}^{(\alpot-1 , l+1)}(2r^{2} - 1) \nonumber \\
 & \mbox{ }  \quad  \quad   \ - \
( \mp ) (n + 1) (1 - r^{2})^{\alpot - 1} \mcV_{l-1 , \mu^{*}}(\bfx) P_{n+1}^{(\alpot-1 , l-1)}(2r^{2} - 1) \, .    \label{ders2}  
 \end{align}
 \end{theorem}

 \begin{theorem} \cite{zhe232} \label{figradinR2}
 Let $f(\bfx) \ = \ (1 - r^{2})^{\frac{\alpha}{2}} \mcV_{l , \mu}(\bfx) P_{n}^{(\frac{\alpha}{2} , l)}(2r^{2} - 1)$. Then, 
 \begin{align*}
(- \Delta)^{\frac{\alpha - 2}{2}} \, \frac{\partial f}{\partial x} &= \
  C_{2} \, (n + \alpot + l) \, \mcV_{l+1 , \mu}(\bfx) P_{n}^{(\alpot-1 , l+1)}(2r^{2} - 1) \nonumber \\
 & \mbox{ } \quad \quad \ + \
C_{2} \, (n + l + 1) \, \mcV_{l-1 , \mu}(\bfx) P_{n+1}^{(\alpot-1 , l-1)}(2r^{2} - 1) \, ,   \\   
(- \Delta)^{\frac{\alpha - 2}{2}} \, \frac{\partial f}{\partial y} &= \
  C_{2} \,  ( \pm ) (n + \alpot + l) \, \mcV_{l+1 , \mu^{*}}(\bfx) P_{n}^{(\alpot-1 , l+1)}(2r^{2} - 1)    \\   
 & \mbox{ }  \quad  \quad   \ + \
 C_{2} \, ( \mp ) (n + l + 1) \, \mcV_{l-1 , \mu^{*}}(\bfx) P_{n+1}^{(\alpot-1 , l-1)}(2r^{2} - 1) \, ,   \\  
\mbox{ where } C_{2} &= \  - 2^{\alpha - 2} \,   \frac{ \Gamma(n + 1 + \frac{\alpha}{2}) \, \Gamma(n + \frac{\alpha}{2} + l)}%
  {\Gamma(n + 1)\, \Gamma(n + 2 + l)}    \, .                    
 \end{align*}
 \end{theorem}

\begin{corollary} \cite{zhe232} \label{gradwR2}
 Let $f(\bfx) \ = \ \mcV_{l , \mu}(\bfx) P_{n}^{(\gamma , l)}(2r^{2} - 1)$. Then,
 \begin{align}
\frac{\partial f}{\partial x} &= \
  (n + l)  \, \mcV_{l-1 , \mu}(\bfx) P_{n}^{(\gamma + 1 \, , \, l - 1)}(2r^{2} - 1)   \nonumber \\
& \hspace{2.0in} 
 \ + \
(n + \gamma + l + 1) \,  \mcV_{l+1 , \mu}(\bfx) P_{n-1}^{(\gamma + 1 \, , \, l + 1)}(2r^{2} - 1) \, ,    \label{dersw1} \\
\mbox{and} &  \nonumber \\
\frac{\partial f}{\partial y} &= \
 (\mp) (n + l) \,  \mcV_{l-1 \, , \mu^{*}}(\bfx) \, P_{n}^{(\gamma + 1 \, , \,  l - 1)}(2 r^{2} - 1)    \nonumber \\
 & \hspace{2.0in} 
\ + \   (\pm) (n + \gamma + l + 1) \, \mcV_{l+1 \, , \mu^{*}}(\bfx) \,   P_{n-1}^{(\gamma + 1 \,  ,  \, l+1)}(2 r^{2} - 1)  \, .  
  \label{dersw2}  
 \end{align}
 \end{corollary}

\subsubsection{Action of the fractional Laplacian}
\label{sssec_fLap}   

We have the following result for the fractional Laplacian.
\begin{theorem} \cite[Theorem 3]{dyd171}  \label{genThmLapm}
For $\alpha > 0$, $l, n \, \geq 0$  integer, $\bfx \in \Omega$, 
\begin{align}
   \left( -\Delta \right)^{\frac{\alpha}{2}} (1 \, - \, r^{2})^{\alpot} \, \mcV_{l , \mu}(\bfx) \, P_{n}^{(\alpot ,  l)}(2 r^{2} \, - \, 1) &= \   
   \lambda_{l, n} \, 
  \mcV_{l , \mu}(\bfx) \, P_{n}^{(\alpot ,  l)}(2 r^{2} \, - \, 1) \, ,    \label{bceq0v}  \\
\mbox{where } \ \lambda_{l, n} &= \   2^{\alpha} \, 
  \frac{ \Gamma(n + 1 + \frac{\alpha}{2}) \, \Gamma(n + l + 1 + \frac{\alpha}{2})}%
  {\Gamma(n + 1)\, \Gamma(n + l + 1)} \, .    \nonumber 
\end{align}
\end{theorem}

\subsubsection{Action of the Riesz potential operator}
\label{sssec_RPO}     

The following theorem, used in the analysis below, is an extension of Theorem \ref{genThmLapm}
from the fractional Laplace operator to the
Riesz potential operator.
\begin{theorem} \cite{zhe232}  \label{genThm3v2}
For $\alpha > 0$, $l, n \, \geq 0$  integer, $s$ an integer, $\alpot - s > -1$, $\bfx \in \Omega$, and   \linebreak[4]
$f(\bfx) \ = \ (1 \, - \, r^{2})^{\alpot - s}_{+} \, \mcV_{l , \mu}(\bfx) \, P_{n}^{(\alpot - s \, , \,  l)}(2 r^{2} \, - \, 1)$,
\[
   \left( -\Delta \right)^{\frac{\alpha - 2}{2}} f(\bfx) \ = \ (-1)^{1 - s} \,  2^{\alpha - 2} \, 
  \frac{ \Gamma(n + 1 - s + \frac{\alpha}{2}) \, \Gamma(n + l + \frac{\alpha}{2})}%
  {\Gamma(n + 1)\, \Gamma(n + l + 2 - s)} \, 
  \mcV_{l , \mu}(\bfx) \, P_{n + 1 - s}^{(\alpot - 2 + s \, ,  l)}(2 r^{2} \, - \, 1) \, .
\]
\end{theorem}   
  
For notation convenience,  we let $\mathbb{N}_{0}  := \mathbb{N} \cup {0}$ and
 use $a \sim b$ to denote that there exists constants $C_{0}$ and $C_{1} > 0$ such that 
 $C_{0} \, b \, \le \, a \,  \le \, C_{1} \, b$. Additionally, we use $a \, \lesssim \, b$ to denote that there exists a constant $C_{1} > 0$ such that
 $a \, \le \, C_{1} \,  b$. 
 
 For $s \in \mathbb{R}$, $\floor{s}$ is used to denote the largest integer that is less than or equal to $s$, and
 $\ceil{s}$ is used to denote the smallest integer that is greater than or equal to $s$.
 
From Stirling's formula we have that
\begin{equation}
\lim_{n \rightarrow \infty} \, \frac{\Gamma(n + \sigma)}{\Gamma(n) \, n^{\sigma}}
\ = \ 1 \, , \mbox{ for } \sigma \in \mathbb{R}.  
 \label{eqStrf}
\end{equation} 
 

 \setcounter{equation}{0}
\setcounter{figure}{0}
\setcounter{table}{0}
\setcounter{theorem}{0}
\setcounter{lemma}{0}
\setcounter{corollary}{0}
\setcounter{definition}{0}
\section{Mapping properties of $(-\Delta)^{\frac{\alpha - 2}{2}}$, and $\Grad$}
\label{sec_map}
In this section we present several mapping properties of the operators. The first mapping property is for the 
Riesz potential operator.

\begin{corollary} \label{cormapD1}
For $\alpha > 0$, we have the following mapping property for $(-\Delta)^{\frac{\alpha - 2}{2}}$ : 
\[
  (-\Delta)^{\frac{\alpha - 2}{2}} \, : \, \omega^{\alpot - 1} \otimes H^{t}_{\alpot - 1}(\Omega) \longrightarrow 
    H^{t + 2 - \alpha}_{\alpot - 1}(\Omega)   \, .
\]
Additionally the mapping is bounded, onto and has a bounded inverse.
\end{corollary}
\textbf{Proof}: Let $\wtilde{v}(\bfx) \ =  \ \omega^{\alpot - 1} v(\bfx)$, where
\be
v(\bfx) \ = \ \sum_{l, n, \mu} a_{l, n, \mu} \, \mcV_{l, \mu}(\bfx) P_{n}^{(\alpot - 1 \,  , l)}(\rho) \ \in H^{t}_{\alpot - 1}(\Omega) \, ,
\label{ert1}
\ee
i.e.,
\[
  \sum_{l, n, \mu} (n + 1)^{t} \, (n + l + 1)^{t} \, 
  a_{l, n, \mu}^{2} \, \| \mcV_{l, \mu}(\bfx) P_{n}^{(\alpot - 1 \,  , l)}(\rho) \|_{L^{2}_{\alpot - 1}}^{2} \ < \ \infty \, .
\]

From Theorem \ref{genThm3v2} (with $s = 1$) we have
\[
 (-\Delta)^{\frac{\alpha - 2}{2}} \, \omega^{\alpot - 1} \, \mcV_{l, \mu}(\bfx) P_{n}^{(\alpot - 1 \,  , l)}(\rho)
 \ = \ 2^{\alpha - 2} \, 
  \frac{ \Gamma(n  + \frac{\alpha}{2}) \, \Gamma(n + l + \frac{\alpha}{2})}%
  {\Gamma(n + 1)\, \Gamma(n + l + 1)} \, 
  \mcV_{l , \mu}(\bfx) \, P_{n}^{(\alpot - 1 \, ,  l)}(\rho) \, .
 \]
 Hence, for $w(\bfx) \, := \,  (-\Delta)^{\frac{\alpha - 2}{2}} \, \omega^{\alpot - 1} \, v(\bfx)$ we have
 \be
 w(\bfx) \ = \ \sum_{l, n, \mu} \, a_{l, n, \mu} \, 2^{\alpha - 2} \, 
  \frac{ \Gamma(n  + \frac{\alpha}{2}) \, \Gamma(n + l + \frac{\alpha}{2})}%
  {\Gamma(n + 1)\, \Gamma(n + l + 1)} \, 
  \mcV_{l , \mu}(\bfx) \, P_{n}^{(\alpot - 1 \, ,  l)}(\rho) \, ,
\label{ert2}
\ee
and (using Stirling's formula \eqref{eqStrf})
\begin{align}
\| w \|_{H^{t + 2 - \alpha}_{\alpot - 1}}^{2} 
&= \  \sum_{l, n, \mu} (n + 1)^{t + 2 - \alpha} \, (n + l + 1)^{t + 2 - \alpha} \,
  \left( 2^{\alpha - 2} \, 
  \frac{ \Gamma(n  + \frac{\alpha}{2}) \, \Gamma(n + l + \frac{\alpha}{2})}%
  {\Gamma(n + 1)\, \Gamma(n + l + 1)} \, a_{l, n, \mu} \right)^{2}        \nonumber \\
 & \hspace{2.0in} \cdot 
   \| \mcV_{l, \mu}(\bfx) P_{n}^{(\alpot - 1 \,  , l)}(\rho) \|_{L^{2}_{\alpot - 1}}^{2}     \nonumber \\
&\sim \    \sum_{l, n, \mu} (n + 1)^{t + 2 - \alpha} \, (n + l + 1)^{t + 2 - \alpha} \, 
\left( (n + 1)^{\alpot - 1} \, (n + l + 1)^{\alpot - 1}  \, a_{l, n, \mu} \right)^{2}       \nonumber \\
 & \hspace{2.0in} \cdot 
   \| \mcV_{l, \mu}(\bfx) P_{n}^{(\alpot - 1 \,  , l)}(\rho) \|_{L^{2}_{\alpot - 1}}^{2}     \nonumber \\
&= \  \sum_{l, n, \mu} (n + 1)^{t} \, (n + l + 1)^{t } \, a_{l, n, \mu}^{2} \, 
\| \mcV_{l, \mu}(\bfx) P_{n}^{(\alpot - 1 \,  , l)}(\rho) \|_{L^{2}_{\alpot - 1}}^{2}    \nonumber \\
&= \ \| v \|_{H^{t}_{\alpot - 1}}^{2} \, .  \label{ert3}
\end{align}  
The explicit relationship between $w(\cdot) \in  H^{t + 2 - \alpha}_{\alpot - 1}(\Omega)$ and
$v(\cdot) \in H^{t}_{\alpot - 1}(\Omega)$ (given by \eqref{ert2} and \eqref{ert1}) establishes the 
onto property of the mapping, and the existence of the inverse. The boundedness of the mapping
and its inverse follow from \eqref{ert3}. \\
\mbox{  } \hfill \qed

The next two results address the mapping properties of the $\Grad$ operator.
\begin{corollary} \label{cormapD2}
We have the following mapping property for the $\Grad$ operator: 
\[
  \Grad \, : \,  H^{s}_{\gamma}(\Omega) \longrightarrow 
   \left( H^{s-1}_{\gamma+1}(\Omega)  \right)^{2} \, .
\]
Additionally the mapping is bounded.
\end{corollary}
\textbf{Remark}: The proof uses Corollary \ref{gradwR2}, and is done in a similar manner as the proof 
of Corollary \ref{cormapD1}. It is given in Appendix \ref{sec_mapPrf}.

The following corollary differs from the previous with respect to the spaces involved.
\begin{corollary} \label{cormapG}
We have the following mapping property for the $\Grad$ operator: 
\[
  \Grad \, : \, \omega^{\gamma} \otimes H^{s}_{\gamma}(\Omega) \longrightarrow 
  \omega^{\gamma-1} \otimes \left( H^{s-1}_{\gamma-1}(\Omega)  \right)^{2} \, .
\]
Additionally the mapping is bounded, one-to-one, with a bounded inverse from its image space.
\end{corollary}
\textbf{Remark}: The proof uses Corollary \ref{gradinR2}, and is done in a similar manner as the proof 
of Corollary \ref{cormapD1}. It is given in Appendix \ref{sec_mapPrf}.

Combining Corollaries \ref{cormapD1} and \ref{cormapG} we have the following.
\begin{corollary} \label{cormapdivG}
We have the following mapping property for the $(-\Delta)^{\frac{\alpha - 2}{2}} \Grad$ operator: 
\[
  (-\Delta)^{\frac{\alpha - 2}{2}} \Grad \, : \, \omega^{\alpot} \otimes H^{s}_{\alpot}(\Omega) \longrightarrow 
  \left( H^{s+1-\alpha}_{\alpot-1}(\Omega) \right)^{2} \, .
\]
Additionally the mapping is bounded, one-to-one, with a bounded inverse from its image space.
\end{corollary}
\mbox{  } \hfill \qed

 \setcounter{equation}{0}
\setcounter{figure}{0}
\setcounter{table}{0}
\setcounter{theorem}{0}
\setcounter{lemma}{0}
\setcounter{corollary}{0}
\setcounter{definition}{0}
\section{Existence and uniqueness of solution}
\label{sec_exu}
With the family of functions defined in Section \ref{ssec_H01}, and the mapping properties established in
Section \ref{sec_map}, we are now in a position to discuss the existence and uniqueness of solutions to 
\eqref{defLL}.
We seek a solution $\wtilde{u}$ 
with the (leading) singular component explicitly represented in the form of the solution. Hence we assume a solution
of the form $\wtilde{u}(\bfx) \, = \, \omega^{\alpot} u(\bfx)$, with $u(\bfx) \in H^{s}_{\alpot}(\Omega)$. With this in mind
we introduce the following additional notation. Let
\begin{align*}
\mcH^{s}_{\beta}(\Omega) &:= \, \omega^{\beta} \otimes H_{\beta}^{s}(\Omega) \, = \, 
\left\{ \wtilde{v}(\bfx) \, = \, \omega^{\beta} v(\bfx) \, : \, v(\bfx) \in H_{\beta}^{s}(\Omega)  \right\} \, ,  \ \ 
\mbox{with } \ \| \wtilde{v} \|_{\mcH^{s}_{\beta}(\Omega)} \, := \, \| v \|_{H_{\beta}^{s}(\Omega)} \, ,  \\
\mbox{and } 
\big( \mcH^{s}_{\beta}(\Omega) \big)^{2} &:= \, \omega^{\beta} \otimes \big( H_{\beta}^{s}(\Omega) \big)^{2} \, = \, 
\left\{ \wtilde{V}(\bfx) \, = \, \omega^{\beta} [V_{1}(\bfx) , V_{2}(\bfx)]^{T} \, : \, V_{1}(\bfx) , \, V_{1}(\bfx) \in 
H_{\beta}^{s}(\Omega)  \right\} \, ,   \\
& \quad \quad \mbox{with } 
 \| \wtilde{V} \|^{2}_{\big( \mcH^{s}_{\beta}(\Omega) \big)^{2}} \, := \, \| V_{1} \|_{H_{\beta}^{s}(\Omega)}^{2} \, + \, 
 \| V_{2} \|_{H_{\beta}^{s}(\Omega)}^{2}   \, .
 \end{align*}

Multiplying \eqref{defLL} through by $\wtilde{v}\, = \, \omega^{\alpot} v(\bfx)$ and integrating over $\Omega$ (using Lemma \ref{mvRPO})
we obtain
\be
B(\wtilde{u} , \wtilde{v}) \, := \, \left( K \Grad \wtilde{u} \, , \, (-\Delta)^{\frac{\alpha - 2}{2}} \Grad \wtilde{v} \right) \ = \ 
\left(f \, , \, \wtilde{v} \right) \, := \, F(\wtilde{v}) \, .
\label{seq1}
\ee

With $v \in H^{\alpha - 1}_{\alpot}(\Omega)$, 
$\wtilde{u}(\bfx) \, = \, \omega^{\alpot} u(\bfx)$, with $u \in H^{1}_{\alpot}(\Omega)$, and 
$K \in \big( W_{w}^{0, \infty}(\Omega) \big)^{2 \times 2}$,
from the mapping properties given in Corollaries \ref{cormapdivG}, \ref{cormapG} the inner product on the LHS of \eqref{seq1}
is well defined. Additionally, for $f \in H_{\alpot}^{-(\alpha - 1)}(\Omega)$ the inner product on the RHS, interpreted as a duality pairing
is also well defined. 
Because of the different test and trial spaces, the analysis of existence and uniqueness of 
the solution to the weak formulation  \eqref{seq1} 
is most appropriately studied via the Banach-Ne\v{c}as-Babu\v{s}ka Theorem \cite[Theorem 2.6]{ern041}. The bilinear form
$B(\cdot , \cdot)$ must satisfy the following three properties. There exists constants $C_{1}, \, C_{2} > 0$ such that
\begin{align}
 B(\wtilde{u} \, , \, \wtilde{z}) \ \le \ C_{1} & \| \wtilde{u} \|_{ \mcH_{\alpot}^{1}(\Omega)} \, 
  \| \wtilde{z} \|_{ \mcH_{\alpot}^{\alpha - 1}(\Omega)} \, , \ \ \forall  \ \ \wtilde{u} \in \mcH_{\alpot}^{1}(\Omega) \, ,  \ 
     \wtilde{z} \in \mcH_{\alpot}^{\alpha - 1}(\Omega) \, ,   \label{condd1}  \\
   \sup_{ \wtilde{z} \in \mcH_{\alpot}^{\alpha - 1}(\Omega)} 
   \frac{ | B(\wtilde{u} \, , \, \wtilde{z})  |}{ \| \wtilde{z} \|_{ \mcH_{\alpot}^{\alpha - 1}(\Omega)}}  &\ge \ 
     C_{2} \,  \| \wtilde{u} \|_{ \mcH_{\alpot}^{1}(\Omega)} \, , \ \ \forall \ \  \wtilde{u} \in \mcH_{\alpot}^{1}(\Omega) \, ,   \label{condd2}  \\
  \sup_{\wtilde{u} \in \mcH_{\alpot}^{1}(\Omega)} 
    | B(\wtilde{u} \, , \, \wtilde{z})  |   &> \ 
     0 \, , \ \ \forall  \ \ 0 \neq \wtilde{z} \in \mcH_{\alpot}^{\alpha - 1}(\Omega) \, .  \label{condd3}
\end{align}

Consistent with the usual variable coefficient diffusion operator, we will assume the $K(\cdot)$ is a symmetric positive definite
matrix satisfying, for $\lambda_{m} , \, \lambda_{M} \in \real^{+}$,
\be
\lambda_{m} \, \xi^{T} \, \xi \ \le \ \xi^{T} \, K(\bfx) \, \xi \ \le \lambda_{M} \, \xi^{T} \, \xi \ \le \ \xi^{T} \, , 
\mbox{ for all } \xi \in \real^{2}, \ \bfx \in \Omega \, .
\label{Kspd}
\ee

The following three lemmas are used in establishing the existence and uniqueness of the solution to \eqref{seq1}.
\begin{lemma} \label{lma_nmeq1}
Let $\wtilde{u} \in \mcH^{1}_{\alpot}(\Omega)$ and 
$U \in \big( H^{0}_{\alpot - 1}(\Omega) \big)^{2} \, = \, \big( L^{2}_{\alpot - 1}(\Omega) \big)^{2}$ be given by
$\Grad \wtilde{u} \, = \, \omega^{\alpot - 1} \, U$. Then
\be 
  \frac{9  \alpha \, + \, 2}{8}  \, \| \wtilde{u} \|^{2}_{\mcH^{1}_{\alpot}(\Omega)} \ \le \ \| U \|^{2}_{\big( L^{2}_{\alpot - 1}(\Omega) \big)^{2}}
     \ \le \  4 \,  \| \wtilde{u} \|^{2}_{\mcH^{1}_{\alpot}(\Omega)} \, .
\label{Hnmeq1} 
\ee
\end{lemma}
\textbf{Proof}: 
\begin{align}
\mbox{Let } \ \ \wtilde{u}(\bfx) &= \ \omega^{\alpot} \, \sum_{l, n, \mu} a_{l, n, \mu} \, \mcV_{l, \mu}(\bfx) \, P_{n}^{(\alpot , l)}(\rho) \, .
   \label{hnb0} \\
  \mbox{Then  } \ \ 
 \| \wtilde{u} \|^{2}_{\mcH^{1}_{\alpot}(\Omega)} &= \ 
    \sum_{l, n, \mu} (n + 1) \, (n + l + 1) \, a^{2}_{l, n, \mu} \, \| \mcV_{l, \mu}(\bfx) \, P_{n}^{(\alpot , l)}(\rho) \|^{2}_{L^{2}_{\alpot}}   \, .
  \label{hnb1}
\end{align}  

Using Theorem \ref{gradinR2}, keeping in mind that $a_{l, n, \mu} = 0$ for $l, n < 0$, $a_{0, n, -1} = 0$, 
$\mcV_{0, -1}(\bfx) = 0$, after reindexing, 
\begin{align}
U_{1} &= \ \sum_{l, n, \mu} - \Big( (n + \alpot) \, a_{l-1, n, \mu} \ + \ n \, a_{l+1, n-1, \mu} \Big) \, 
    \mcV_{l, \mu}(\bfx) \, P_{n}^{(\alpot - 1 \,  , \, l)}(\rho)  \, ,    \label{defU1} \\
\mbox{and} \     
U_{2} &= \ \sum_{l, n, \mu} - \Big( (\pm) (n + \alpot) \, a_{l-1, n, \mu} \ + \ (\mp) n \, a_{l+1, n-1, \mu} \Big) \, 
    \mcV_{l, \mu^{*}}(\bfx) \, P_{n}^{(\alpot - 1 \,  , \, l)}(\rho)  \,  . \label{defU2}
\end{align}

Using orthogonality of the basis functions $\mcV_{l, \mu}(\bfx) \, P_{n}^{(\alpot - 1 \,  , \, l)}(\rho)$ in $L^{2}_{\alpot - 1}(\Omega)$, 
\linebreak[4] 
$\| \mcV_{l, \mu}(\bfx) \, P_{n}^{(\alpot - 1 \,  , \, l)}(\rho) \|_{L^{2}_{\alpot - 1}(\Omega)} \ = \ 
  \| \mcV_{l, \mu^{*}}(\bfx) \, P_{n}^{(\alpot - 1 \,  , \, l)}(\rho) \|_{L^{2}_{\alpot - 1}(\Omega)}$, and that
 formally $\mcV_{0, -1}(\bfx) \, P_{n}^{(\alpot - 1 \,  , \, l)}(\rho)$ are not basis functions, we obtain
 \begin{align}
 & \| U_{1} \|_{L^{2}_{\alpot - 1}(\Omega)}^{2} \, + \, 
 \| U_{2} \|_{L^{2}_{\alpot - 1}(\Omega)}^{2}   \nonumber \\
&= \ \sum_{n \ge 1} n^{2} \, 
\big( a_{1, n-1, 1}^{2} \, + \, a_{1, n-1, -1}^{2} \big) \,  \| \mcV_{0, 1}(\bfx) \, 
P_{n}^{(\alpot - 1 \,  , \, 0)}(\rho) \|^{2}_{L^{2}_{\alpot - 1}} \nonumber \\
& \quad  \ + \ 
 \sum_{l \ge 1, \, n, \mu} 2 \Big( (n + \alpot)^{2} \, a_{l-1, n, \mu}^{2} \ + \ n^{2} \, a_{l+1, n-1, \mu}^{2} \, \Big) \,  
  \| \mcV_{l, \mu}(\bfx) \, P_{n}^{(\alpot - 1 \,  , \, l)}(\rho) \|^{2}_{L^{2}_{\alpot - 1}}  \nonumber \\
&= \ \sum_{n \ge 1} n^{2} \, 
\big( a_{1, n-1, 1}^{2}  \, + \, a_{1, n-1, -1}^{2} \big) \,  \| \mcV_{0, 1}(\bfx) \, 
P_{n}^{(\alpot - 1 \,  , \, 0)}(\rho) \|^{2}_{L^{2}_{\alpot - 1}}  \nonumber \\
& \ + \ 
2  \sum_{l \ge 1, \, n, \mu}  (n + \alpot)^{2} \, a_{l-1, n, \mu}^{2}  \,  
  \| \mcV_{l, \mu}(\bfx) \, P_{n}^{(\alpot - 1 \,  , \, l)}(\rho) \|^{2}_{L^{2}_{\alpot - 1}}     \nonumber \\
& \ + \
2  \sum_{l \ge 1, \, n \ge 1, \mu}  n^{2} \, a_{l+1, n-1, \mu}^{2} \,  
  \| \mcV_{l, \mu}(\bfx) \, P_{n}^{(\alpot - 1 \,  , \, l)}(\rho) \|^{2}_{L^{2}_{\alpot - 1}}   \, .  \label{hnb2} 
 \end{align}

Using Lemma \ref{lmaeq5} and Lemma \ref{lmaeq3},
\be
\| \mcV_{l, \mu}(\bfx) \, P_{n}^{(\alpot - 1 \,  , \, l)}(\rho) \|^{2}_{L^{2}_{\alpot - 1}} \ = \ 
\frac{C_{l, \mu}}{C_{l+1, \mu}} \, \frac{(n + l + \alpot)}{n} \,
 \| \mcV_{l+1, \mu}(\bfx) \, P_{n-1}^{(\alpot  , \, l+1)}(\rho) \|^{2}_{L^{2}_{\alpot}} \, ,
 \label{hnb3}
 \ee
 and using Lemma \ref{lmaeq4}  and Lemma \ref{lmaeq3},
\be
\| \mcV_{l, \mu}(\bfx) \, P_{n}^{(\alpot - 1 \,  , \, l)}(\rho) \|^{2}_{L^{2}_{\alpot - 1}} \ = \ 
\frac{C_{l, \mu}}{C_{l-1, \mu}} \, \frac{(n + l)}{(n + \alpot)} \,
 \| \mcV_{l-1, \mu}(\bfx) \, P_{n}^{(\alpot , \, l-1)}(\rho) \|^{2}_{L^{2}_{\alpot}} \, .
 \label{hnb4}
 \ee
 
 Substituting \eqref{hnb3} and \eqref{hnb4} into \eqref{hnb2},
 \begin{align}
 & \| U_{1} \|_{L^{2}_{\alpot - 1}(\Omega)}^{2} \, + \, 
 \| U_{2} \|_{L^{2}_{\alpot - 1}(\Omega)}^{2}   \nonumber \\
&= \ \sum_{n \ge 1} n^{2} \,\big( a_{1, n-1, 1}^{2} + \, a_{1, n-1, -1}^{2} \big) \,  
\frac{C_{0, 1}}{C_{1, 1}} \, \frac{(n + \alpot)}{n} \, \| \mcV_{1, 1}(\bfx) \, P_{n-1}^{(\alpot \,  , \, 1)}(\rho) \|^{2}_{L^{2}_{\alpot }}  \nonumber \\
& \ + \ 
2  \sum_{l \ge 1, \, n, \mu}  (n + \alpot)^{2} \, 
 a_{l-1, n, \mu}^{2}  \,   \,  \frac{C_{l, \mu}}{C_{l-1 , \mu}} \, \frac{(n + l)}{(n + \alpot)} \, 
  \| \mcV_{l-1 , \mu}(\bfx) \, P_{n}^{(\alpot  \,  , \, l-1)}(\rho) \|^{2}_{L^{2}_{\alpot}}     \nonumber \\
& \ + \
2  \sum_{l \ge 1, \, n \ge 1, \mu}  n^{2} \, a_{l+1, n-1, \mu}^{2} \,  \,  \frac{C_{l, \mu}}{C_{l+1 , \mu}} \, \frac{(n + l + \alpot)}{n} \, 
  \| \mcV_{l+1, \mu}(\bfx) \, P_{n-1}^{(\alpot \,  , \, l+1)}(\rho) \|^{2}_{L^{2}_{\alpot}}   \nonumber  \\
&= \ \sum_{n \ge 0} (n + 1) \, (n + \alpot + 1) \, \frac{1}{2} \, 
\big( a_{1, n, 1}^{2}  \, + \, a_{1, n, -1}^{2} \big) \,  
 \| \mcV_{1, 1}(\bfx) \, P_{n}^{(\alpot \,  , \, 1)}(\rho) \|^{2}_{L^{2}_{\alpot }}  \nonumber \\
& \ + \ 
2  \sum_{n \ge 0}  (n + \alpot) \, (n + 1) \, 2 \, a_{0, n, 1}^{2}  \,  
  \| \mcV_{0 , 1}(\bfx) \, P_{n}^{(\alpot  \,  , \, 0)}(\rho) \|^{2}_{L^{2}_{\alpot}}     \nonumber \\
& \ + \ 
2  \sum_{l \ge 1, \, n, \mu}  (n + \alpot) \, (n + l + 1) \, a_{l, n, \mu}^{2}  \, 
  \| \mcV_{l , \mu}(\bfx) \, P_{n}^{(\alpot  \,  , \, l)}(\rho) \|^{2}_{L^{2}_{\alpot}}     \nonumber \\
& \ + \
2  \sum_{l \ge 2, \, n \ge 0, \mu}  (n + 1) \, (n + l + \alpot) \, a_{l, n, \mu}^{2}  \, 
  \| \mcV_{l, \mu}(\bfx) \, P_{n}^{(\alpot \,  , \, l)}(\rho) \|^{2}_{L^{2}_{\alpot}}    \label{hnb5}
 \end{align}
 
Combining the appropriate terms in the various summations in \eqref{hnb5} we obtain
\begin{align}
& \| U_{1} \|_{L^{2}_{\alpot - 1}(\Omega)}^{2} \, + \, 
 \| U_{2} \|_{L^{2}_{\alpot - 1}(\Omega)}^{2}   \nonumber \\
&= \  \sum_{l = 0 , \, n \ge 0 \, , \mu = 1}  4 \, (n + \alpot) \, (n + 1) \, a_{0, n, 1}^{2}  \,  
  \| \mcV_{0 , 1}(\bfx) \, P_{n}^{(\alpot  \,  , \, 0)}(\rho) \|^{2}_{L^{2}_{\alpot}}     \nonumber \\
& \ + \
\sum_{l = 1 , \, n \ge 0 \, , \mu} \Big( \frac{1}{2} \, (n + 1) \, (n + \alpot + 1) \ + \
   2 \, (n + \alpot) \, (n + 2) \Big) \,  a_{1, n, \mu}^{2} \,  
 \| \mcV_{1, \mu}(\bfx) \, P_{n}^{(\alpot \,  , \, 1)}(\rho) \|^{2}_{L^{2}_{\alpot }}  \nonumber \\
& \ + \
 \sum_{l \ge 2, \, n \ge 0, \mu} 2 \, \Big(  (n + \alpot) \, (n + l + 1)  \ + \ (n + 1) \, (n + l + \alpot) \Big) \, a_{l, n, \mu}^{2}  \, 
  \| \mcV_{l, \mu}(\bfx) \, P_{n}^{(\alpot \,  , \, l)}(\rho) \|^{2}_{L^{2}_{\alpot}}   \, . \label{hnb6}
\end{align}

Comparing the coefficients in \eqref{hnb6} with \eqref{hnb1} and noting that
\begin{align}
  2 \alpha &\le \ \frac{4 \, (n + \alpot) \, (n + 1)}{(n + 1) \, (n + 1)} \ < \ 4 \, , \ \mbox{ for } n \ge 0 \, ,   \label{rato1} \\
\frac{9 \alpha \, + \, 2}{8} &\le  \ \frac{ \frac{1}{2} (n + 1) \, (n + \alpot + 1) \ + \ 2 (n + \alpot) \, (n + 2)}{(n + 1) \, (n + 2)} \ 
< \ \frac{5}{2} \, , \ \mbox{ for } n \ge 0 \, ,   \label{rato2} \\
%
 \Big( \frac{4 \alpha \, + \, 4 }{3} \Big) &\le \ \frac{ 2 \big( (n + \alpot) \, (n + l + 1) \ + \ (n + 1) \, (n + l + \alpot) \big)}%
 {(n + 1) \, (n + l + 1)} \ < \ 4  \, , \ \mbox{ for } n \ge 0 \, ,  l \ge 2 \, ,   \label{rato4}
\end{align}
we obtain \eqref{Hnmeq1}. \\
\mbox{  } \hfill \qed

\begin{lemma} \label{lmau2v}
Let $\wtilde{u}(\bfx) \ = \ \omega^{\alpot} \, \sum_{l, n, \mu} a_{l, n, \mu} \, \mcV_{l, \mu}(\bfx) \, P_{n}^{(\alpot , l)}(\rho) \, \in \, 
\mcH_{\alpot}^{1}(\Omega)$. Then
\begin{align}
& \wtilde{v}(\bfx) \ = \ 2^{2 - \alpha} \, \omega^{\alpot} \, \sum_{l, n, \mu} (n + 1) \, \frac{\Gamma(n + 1)}{\Gamma(n + 1 + \alpot)} \, 
 \frac{\Gamma(n + l + 1)}{\Gamma(n + l + \alpot)} \, 
 a_{l, n, \mu} \, \mcV_{l, \mu}(\bfx) \, P_{n}^{(\alpot , l)}(\rho) \, \in \, \mcH_{\alpot}^{\alpha - 1}(\Omega) \, ,  \label{defvtda} \\
 & \ \mbox{with } \ \| \wtilde{u} \|_{\mcH_{\alpot}^{1}(\Omega)} \ \sim  \ \| \wtilde{v} \|_{\mcH_{\alpot}^{\alpha - 1}(\Omega)}  \, . \nonumber 
\end{align}
\end{lemma}
\textbf{Proof}: We have
\begin{align*}
& \| \wtilde{v} \|^{2}_{\mcH_{\alpot}^{\alpha - 1}}  \ =   \nonumber \\
& \ 4^{2 - \alpha}
\sum_{l, n, \mu}  (n + 1)^{\alpha - 1} \, (n + l + 1)^{\alpha - 1} \, 
\Big( (n + 1) \, \frac{\Gamma(n + 1)}{\Gamma(n + 1 + \alpot)} \, 
 \frac{\Gamma(n + l + 1)}{\Gamma(n + l + \alpot)} \, 
 a_{l, n, \mu} \Big)^{2} \, \| \mcV_{l, \mu}(\bfx) \, P_{n}^{(\alpot , l)}(\rho) \|^{2}_{L^{2}_{\alpot}}  \\
&\sim \ 
\sum_{l, n, \mu}  (n + 1)^{\alpha - 1} \, (n + l + 1)^{\alpha - 1} \, 
 (n + 1)^{2} \, (n + 1)^{- \alpot \, 2} \, (n + l + 1)^{- (\alpot - 1) \, 2} \, 
 a_{l, n, \mu}^{2} \, \| \mcV_{l, \mu}(\bfx) \, P_{n}^{(\alpot , l)}(\rho) \|^{2}_{L^{2}_{\alpot}}  \\
& \quad \quad  \quad \quad  \mbox{(using Stirling's formula \eqref{eqStrf})}  \\
&= \ 
\sum_{l, n, \mu}  (n + 1)^{1} \, (n + l + 1)^{1} \, 
a_{l, n, \mu}^{2} \, \| \mcV_{l, \mu}(\bfx) \, P_{n}^{(\alpot , l)}(\rho) \|^{2}_{L^{2}_{\alpot}}   \\
&= \ \| \wtilde{u} \|^{2}_{\mcH_{\alpot}^{1}} \, . 
\end{align*}
\mbox{  } \hfill \qed

\begin{lemma} \label{lmaBdDif}
Let $\wtilde{u}(\bfx)$ and $\wtilde{v}(\bfx)$ be defined as in Lemma \ref{lmau2v}, and $U$ and $V$ given by
$\Grad \wtilde{u} \ = \ \omega^{\alpot - 1} \, U$ and $( - \Delta )^{\frac{2 - \alpha}{2}} \Grad \wtilde{v} \ = \ V$.
Then for $W \ = \ U \, - \, V$ we have
\be
\| W \|_{L^{2}_{\alpot - 1}(\Omega)} \ \le \ \frac{(2 - \alpha)}{\sqrt{\alpha \, (2 + \alpha)}} \, \| U \|_{L^{2}_{\alpot - 1}(\Omega)}     \, . 
\label{Wbd1}
\ee
\end{lemma}
\textbf{Proof}:   
Let $U \ = \  [U_{1} , \, U_{2}]^{T}$, where
$U_{1}, \, U_{2}$ are given by \eqref{defU1} and \eqref{defU2}, respectively. Also, let $V \, = \, [V_{1} , \, V_{2}]^{T}$
where, using 
%
Theorem \ref{figradinR2}, after reindexing, we have
\begin{align}
V_{1} &= \ \sum_{l, n, \mu} - \Big( (n + 1) \, \frac{(n + \alpot + l - 1)}{(n + l)} \, a_{l-1, n, \mu} \ + \ n \, a_{l+1, n-1, \mu} \Big) \, 
    \mcV_{l, \mu}(\bfx) \, P_{n}^{(\alpot - 1 \,  , \, l)}(\rho)  \, ,    \label{defV1} \\
\mbox{and} \     
V_{2} &= \ \sum_{l, n, \mu} - \Big( (\pm) (n + 1) \, \frac{(n + \alpot + l - 1)}{(n + l)}  \, a_{l-1, n, \mu} \ + \ (\mp) n \, a_{l+1, n-1, \mu} \Big) \, 
    \mcV_{l, \mu^{*}}(\bfx) \, P_{n}^{(\alpot - 1 \,  , \, l)}(\rho)  \,  . \label{defV2}
\end{align}

Next, let $W \, = \, [W_{1} , \, W_{2}]^{T} \, := \, U \, - \, V$. Subtracting \eqref{defV1} from \eqref{defU1} and 
 \eqref{defV2} from \eqref{defU2} yields
\begin{align*}
W_{1} &= \ \sum_{l \ge 1, n, \mu} \Big( (n + 1) \, \frac{(n + \alpot + l - 1)}{(n + l)} \, - \, (n + \alpot) \Big)
\, a_{l-1, n, \mu} \, 
    \mcV_{l, \mu}(\bfx) \, P_{n}^{(\alpot - 1 \,  , \, l)}(\rho)  \, ,   \\  
\mbox{and} \     
W_{2} &= \ \sum_{l \ge 1, n, \mu} (\pm) \Big(  (n + 1) \, \frac{(n + \alpot + l - 1)}{(n + l)} \, - \, (n + \alpot) \Big) 
 \, a_{l-1, n, \mu}  \, 
    \mcV_{l, \mu^{*}}(\bfx) \, P_{n}^{(\alpot - 1 \,  , \, l)}(\rho)  \,  .   
\end{align*}

Using orthogonality of the basis functions and \eqref{hnb4} we obtain
\begin{align}
 & \| W_{1} \|_{L^{2}_{\alpot - 1}(\Omega)}^{2} \, + \, 
 \| W_{2} \|_{L^{2}_{\alpot - 1}(\Omega)}^{2}   \nonumber \\
&= \  \sum_{l \ge 1, n, \mu} 2 \, \Big( (n + 1) \, \frac{(n + \alpot + l - 1)}{(n + l)} \ - \ (n + \alpot) \Big)^{2} \, 
a_{l-1, n, \mu}^{2} \, \frac{C_{l, \mu}}{C_{l-1, \mu}} \, \frac{(n + l)}{(n + \alpot)} \, 
\| \mcV_{l-1, \mu}(\bfx) \, P_{n}^{(\alpot , \,  l-1)}(\rho) \|_{L^{2}_{\alpot}}^{2}   \nonumber \\
&= \  \sum_{l \ge 1, n, \mu} 2 \,\Big( (n + 1) \, \frac{(n + \alpot + l)}{(n + l + 1)} \ - \ (n + \alpot) \Big)^{2} \,  \frac{(n + l + 1)}{(n + \alpot)}  \, 
a_{l, n, \mu}^{2} \, 
\| \mcV_{l, \mu}(\bfx) \, P_{n}^{(\alpot , \,  l)}(\rho) \|_{L^{2}_{\alpot}}^{2}  \, ,  \label{hnb11}
\end{align}
where we have used the terms for $l = 1$ in the equation before \eqref{hnb11} are $0$.

As the $a_{l, n , \mu}$ are arbitrary, in order
to bound  $\| W \|_{L^{2}_{\alpot - 1}(\Omega)}$ by $\| U \|_{L^{2}_{\alpot - 1}(\Omega)}$  
we must consider the ratio of their individual coefficients in \eqref{hnb6} and \eqref{hnb11}.

For $l = 1, \, n \ge 0$, of interest is
\begin{align}
f_{1}(n) &:= \ \frac{ 2 \, \Big( (n + 1) \, \frac{(n + \alpot + 1)}{(n + 2)} \ - \ (n + \alpot) \Big)^{2} \, \frac{(n + 2)}{(n + \alpot)} }%
{ \frac{1}{2} \, (n + 1) \, (n + \alpot + 1) \ + \ 2 \, (n + \alpot) \, (n + 2) }    \nonumber \\
&= \frac{(x - 1)^{2}}{\frac{1}{4} \, x \ + 1}  \ := \ \hat{f}_{1}(x) \,  ,  \ \mbox{ where }
 x \, = \, \frac{(n + 1) \, (n + \alpot + 1)}{(n + 2) \, (n + \alpot)} \, .    \label{dfx1}
\end{align}
As $\frac{d x}{dn} < 0$, of interest is $\hat{f}_{1}(x)$ for $1 < x \le \, (1 + \alpot) / \alpha$. As $\frac{d \hat{f}_{1}(x)}{dx} > 0$, for
$1 < x \le \, (1 + \alpot) / \alpha$,
\be
\sup_{n \ge 0} f_{1}(n) \ = \ \hat{f}_{1}((1 + \alpot) / \alpha) \ = \ \frac{2 \, (2 - \alpha)^{2}}{\alpha \, (2 \, + \, 9 \alpha)} \, .
\label{f1mx}
\ee

%
 
For $l \ge 2, \, n \ge 0$, of interest is
\begin{align}
f_{3}(l, n) &:= \ \frac{ 2 \, \Big( (n + 1) \, \frac{(n + \alpot + l)}{(n + l + 1)} \ - \ (n + \alpot) \Big)^{2} \, \frac{(n + l + 1)}{(n + \alpot)} }%
{ 2 \, \big(  (n + \alpot) \, (n + l + 1) \ + \  (n + 1) \, (n + l + \alpot) \big) }    \nonumber \\
&= \frac{(x - 1)^{2}}{(x \, + \,1)}  \ := \ \hat{f}_{3}(x) \,  ,  \ \mbox{ where }
 x \, = \, \frac{(n + 1) \, (n + l + \alpot)}{(n + l + 1) \, (n + \alpot)} \, .    \label{dfx2}
\end{align}
As $\frac{\partial x}{\partial l} > 0$, $\frac{\partial x}{\partial n} < 0$, (for $1 < \alpha < 2$)
of interest is $\hat{f}_{3}(x)$ for $1 < x < \, 2  / \alpha$. As $\frac{d \hat{f}_{3}(x)}{dx} > 0$, for
$1 < x < \, 2  / \alpha $,
\[
\sup_{l \ge 2, \, n \ge 0} f_{3}(l, n) \ = \ \hat{f}_{3}(2  / \alpha) \ = \ \frac{(2 - \alpha)^{2}}{\alpha \, (2 + \alpha)} \, .
\]
%
Noting that
\[
       \sup_{1 < \alpha < 2} \left\{ \frac{2 \, (2 - \alpha)^{2}}{\alpha \, (2 \, + \, 9 \alpha)}  \, , \ 
    \frac{(2 - \alpha)^{2}}{\alpha \, (2 + \alpha)}  \right\}
     \ = \    \frac{(2 - \alpha)^{2}}{\alpha \, (2 + \alpha)}   \, ,
\]     
it follows that
\[
\| W \|_{L^{2}_{\alpot - 1}(\Omega)} \ \le \ \frac{(2 - \alpha)}{\sqrt{\alpha \, (2 + \alpha)}}  \, \| U \|_{L^{2}_{\alpot - 1}(\Omega)}     \, . 
\]
\mbox{  } \hfill \qed

We are now in position to establish the existence and uniqueness of solution to \eqref{seq1}.
\begin{theorem}  \label{thmexun}
Assume that $1 < \alpha < 2$, and $K \in \big( W_{w}^{0, \infty}(\Omega) \big)^{2 \times 2}$ is a symmetric positive definite matrix
with $\lambda_{m}$ and $\lambda_{M}$ in \eqref{Kspd} satisfying 
$\lambda_{M} \, < \, \frac{\sqrt{\alpha \, (2 + \alpha)}}{(2 - \alpha)} \, \lambda_{m}$.
Then, given $f \in H^{-(\alpha - 1)}_{\alpot}(\Omega)$ there exists a unique $\wtilde{u} \in \mcH^{1}_{\alpot}(\Omega)$ satisfying
\eqref{seq1}, with $\| \wtilde{u} \|_{\mcH^{1}_{\alpot}(\Omega)} \ \le \ \frac{1}{C_{2}} \, \| f \|_{H^{-(\alpha - 1)}_{\alpot}(\Omega)}$,
where $C_{2}$ is the $inf-sup$ constant (see \eqref{condd2}).
\end{theorem}
\textbf{Proof}: 
To establish \eqref{condd1}, for $\wtilde{u} \in \mcH^{1}_{\alpot}(\Omega)$ and $\wtilde{z} \in \mcH^{\alpha -1}_{\alpot}(\Omega)$,
let $U$ and $Z$ be given by $\Grad \wtilde{u} \, = \, \omega^{\alpot - 1} U$ and $(-\Delta)^{\frac{2 - \alpha}{2}} \Grad \wtilde{z} \, = \, Z$.
Then,
\begin{align*}
B(\wtilde{u} \, , \, \wtilde{z}) &= \ \int_{\Omega} K \, \Grad \wtilde{u} \, \cdot \, (-\Delta)^{\frac{2 - \alpha}{2}} \Grad \wtilde{z} \, d\Omega \
= \ \int_{\Omega} K \, \omega^{\alpot - 1}  U \, \cdot \, Z \, d\Omega   \\
&\le \ \big( \int_{\Omega} K \, \omega^{\alpot - 1}  U \, \cdot \, U \, d\Omega \big)^{1/2} \
  \big( \int_{\Omega} K \, \omega^{\alpot - 1}  Z \, \cdot \, Z \, d\Omega \big)^{1/2}   \\
&\le \ \| K \|_{W_{w}^{0, \infty}} \, \| U \|_{L^{2}_{\alpot - 1}} \, \| Z \|_{L^{2}_{\alpot - 1}}  \\
&\le \ C \, \| K \|_{W_{w}^{0, \infty}} \, \| \wtilde{u} \|_{\mcH^{1}_{\alpot}} \,  \| \wtilde{z} \|_{\mcH^{\alpha - 1}_{\alpot}} \, ,
\end{align*}
where in the last step we have used Corollary \ref{cormapG} and Corollary \ref{cormapdivG}.

To establish \eqref{condd2}, let $\wtilde{u}$ be an arbitrary function in $\mcH^{1}_{\alpot}(\Omega)$, and an
associated $\wtilde{v} \in \mcH^{\alpha - 1}_{\alpot}(\Omega)$ be given by \eqref{defvtda}. Let $U$, $V$, and $W$ be given by
$\Grad \wtilde{u} \, = \, \omega^{\alpot - 1} U$, $(-\Delta)^{\frac{2 - \alpha}{2}} \Grad \wtilde{v} \, = \, V$, and $W \, = \, U \, - \, V$. Then,
\begin{align*}
 \sup_{ \wtilde{z} \in \mcH_{\alpot}^{\alpha - 1}} 
   \frac{ | B(\wtilde{u} \, , \, \wtilde{z})  |}{ \| \wtilde{z} \|_{ \mcH_{\alpot}^{\alpha - 1}}}  &\ge \
    \frac{ | B(\wtilde{u} \, , \, \wtilde{v})  |}{ \| \wtilde{v} \|_{ \mcH_{\alpot}^{\alpha - 1}}} \
    \sim \ \frac{ | B(\wtilde{u} \, , \, \wtilde{v})  |}{ \| \wtilde{u} \|_{ \mcH_{\alpot}^{1}}}  \ \mbox{ (using Lemma \ref{lmau2v})}  \\
&= \ \frac{1}{ \| \wtilde{u} \|_{ \mcH_{\alpot}^{1}} } \, \int_{\Omega} K \, \omega^{\alpot - 1} \, U \cdot V \, d\Omega   \\
&= \ \frac{1}{ \| \wtilde{u} \|_{ \mcH_{\alpot}^{1}} } \, 
\Big( \int_{\Omega} K \, \omega^{\alpot - 1} \, U \cdot U \, d\Omega  
\ - \ \int_{\Omega} K \, \omega^{\alpot - 1} \, U \cdot W \, d\Omega  \Big)  \\
&\ge \ \frac{1}{ \| \wtilde{u} \|_{ \mcH_{\alpot}^{1}} } \, 
\big( \lambda_{m} \, \| U \|_{L^{2}_{\alpot - 1}}^{2}  \ - \ \lambda_{M} \, \| U \|_{L^{2}_{\alpot - 1}} \, \| W \|_{L^{2}_{\alpot - 1}} \big) \\
&\ge \ \frac{1}{ \| \wtilde{u} \|_{ \mcH_{\alpot}^{1}} } \, 
\big( \lambda_{m} \ - \ \frac{(2 - \alpha)}{\sqrt{\alpha \, (2 + \alpha)}} \, \lambda_{M} \big)
\, \| U \|_{L^{2}_{\alpot - 1}}^{2}  \ \mbox{ (using Lemma \ref{lmaBdDif})} \\
&\ge \ \frac{9 \alpha \, + \, 2}{8} \, \big( \lambda_{m} \ - \ \frac{(2 - \alpha)}{\sqrt{\alpha \, (2 + \alpha)}} \, \lambda_{M} \big) \, 
\| \wtilde{u} \|_{ \mcH_{\alpot}^{1}}
 \ \mbox{ (using \eqref{Hnmeq1})}  \\
&= \ C_{2} \,  \| \wtilde{u} \|_{ \mcH_{\alpot}^{1}} \, , \ \mbox{ provided } 
\lambda_{M} \, < \, \frac{\sqrt{\alpha \, (2 + \alpha)}}{(2 - \alpha)} \, \lambda_{m} \, .   
\end{align*}

The establishment of \eqref{condd3} is done in a similar manner to \eqref{condd2}.

As \eqref{condd1} - \eqref{condd3} are satisfied, the existence and uniqueness of 
$\wtilde{u} \in \mcH^{1}_{\alpot}(\Omega)$ satisfying \eqref{seq1} follows.

Finally, using \eqref{condd2},
\[
C_{2} \, \| \wtilde{u} \|_{ \mcH_{\alpot}^{1}}  \ \le \  \sup_{ \wtilde{z} \in \mcH_{\alpot}^{\alpha - 1}} 
   \frac{ | B(\wtilde{u} \, , \, \wtilde{z})  |}{ \| \wtilde{z} \|_{ \mcH_{\alpot}^{\alpha - 1}}} 
   \ = \ \sup_{ \wtilde{z} \in \mcH_{\alpot}^{\alpha - 1}} 
   \frac{| ( f , \, \wtilde{z})  |}{ \| \wtilde{z} \|_{ \mcH_{\alpot}^{\alpha - 1}}}    
\ \le \ \frac{ \| f \|_{ \mcH_{\alpot}^{-(\alpha - 1)}}  \, \| \wtilde{z} \|_{ \mcH_{\alpot}^{\alpha - 1}} }{ \| \wtilde{z} \|_{ \mcH_{\alpot}^{\alpha - 1}} }
\ = \ \| f \|_{ \mcH_{\alpot}^{-(\alpha - 1)}}  \, .
\]
\mbox{  } \hfill \qed

\subsection{Higher regularity for $\wtilde{u}$ satisfying \eqref{seq1}}
\label{ssec_Hreg}
With improved regularity of $f$, $f \in H^{- (\alpha - 1) \, + \, s}_{\alpot}(\Omega), \ s \ge 0$, the corresponding solution
of \eqref{seq1} has higher regularity, $\wtilde{u} \in \mcH_{\alpot}^{1 + s}(\Omega)$. For the analysis of this situation,
\eqref{seq1} is replaced by: \\
\textit{Given} $f \in H^{- (\alpha - 1) \, + \, s}_{\alpot}(\Omega), \ s > -1$, 
$K \in \big( W_{w}^{\ceil{s}, \infty}(\Omega) \big)^{2 \times 2}$, symmetric, positive definite, \textit{determine} 
$\wtilde{u} \in \mcH_{\alpot}^{1 + s}(\Omega)$, \textit{such that for all} 
$\wtilde{v} \in \mcH_{\alpot}^{(\alpha - 1) \, - \, s}(\Omega)$
\be
B(\wtilde{u} \, , \, \wtilde{v}) \ := \ 
\langle K \Grad \wtilde{u} \, , \, (-\Delta)^{\frac{\alpha - 2}{2}} \Grad \wtilde{v} \rangle_{[H^{s}_{\alpot} \, , \, H^{-s}_{\alpot}]}
\ = \  \langle f  \, , \,  \wtilde{v} \rangle_{[H^{- (\alpha - 1) \, + \, s}_{\alpot} \, , \, H^{(\alpha - 1) \, - \, s}_{\alpot}]} \, , 
\label{sseq1}
\ee
where $\langle \cdot  \, , \,  \cdot \rangle_{[H^{t}_{\alpot} \, , \, H^{-t}_{\alpot}]}$ denotes the $L^{2}_{\alpot}(\Omega)$
duality pairing between the spaces $H^{t}_{\alpot}(\Omega)$ and $H^{-t}_{\alpot}(\Omega)$.

The proof of the regularity lift for $\wtilde{u}$ follows the same steps as used to establish Theorem \ref{thmexun}.

To simplify the expressions used, we introduce a slight abuse of notation. For $\wtilde{u}$, $\wtilde{v}$, $U$,
$V$, and $W$, as defined in the previous section, let $(n + 1)^{\gamma} \, (n + l + 1)^{\epsilon} \, p(\bfx)$ denote the case
where all the coefficients of $p(\bfx) \in \{ \wtilde{u} , \, \wtilde{v} , \, U , \, V , \, W \}$ relative to the basis for
$L^{2}_{\alpot}(\Omega)$ are multiplied by $(n + 1)^{\gamma} \, (n + l + 1)^{\epsilon}$. For example, 
with $\wtilde{u}(\bfx)$ given by \eqref{hnb0},
\[
(n + 1)^{\gamma} \, (n + l + 1)^{\epsilon} \,  \wtilde{u}(\bfx) 
\ := \ 
\omega^{\alpot} \, \sum_{l, n, \mu} (n + 1)^{\gamma} \, (n + l + 1)^{\epsilon} \, a_{l, n, \mu} \, \mcV_{l, \mu}(\bfx) \, P_{n}^{(\alpot , l)}(\rho) \, .
\]

An important tool in establishing the results in this section is to use the orthogonal basis to redistribute the coefficients
of the basis functions between the product of two functions. We illustrate this tool by the following.
Suppose
$p(\bfx) \ = \ \sum_{l, n, \mu}  \, p_{l, n, \mu} \, \mcV_{l, \mu}(\bfx) \, P_{n}^{(\beta , l)}(\rho)$ and
$q(\bfx) \ = \ \sum_{l, n, \mu}  \, q_{l, n, \mu} \, \mcV_{l, \mu}(\bfx) \, P_{n}^{(\beta , l)}(\rho)$. Then, using orthogonality of the
basis functions,
\begin{align*}
& \ \big( (n + 1)^{\gamma} \, (n + l + 1)^{\epsilon} \, p(\bfx) \, , \, q(\bfx) \big)_{L^{2}_{\beta}}   \\
& \quad = \ 
\int_{\Omega} \omega^{\beta} \, 
\Big( \sum_{l, n, \mu}  \, (n + 1)^{\gamma} \, (n + l + 1)^{\epsilon} \, p_{l, n, \mu} \, \mcV_{l, \mu}(\bfx) \, P_{n}^{(\beta , l)}(\rho) \Big) \
\Big( \sum_{l, n, \mu}  \, q_{l, n, \mu} \, \mcV_{l, \mu}(\bfx) \, P_{n}^{(\beta , l)}(\rho) \Big) \, d\Omega \\
& \quad = \ 
\int_{\Omega} \omega^{\beta} \, 
\Big( \sum_{l, n, \mu}  \, (n + 1)^{\gamma - s} \, (n + l + 1)^{\epsilon - t} \, p_{l, n, \mu} \, \mcV_{l, \mu}(\bfx) \, P_{n}^{(\beta , l)}(\rho) \Big) \
 \cdot \\
& \hspace*{3.0in}
 \Big( \sum_{l, n, \mu}  \, (n + 1)^{s} \, (n + l + 1)^{t}  \, q_{l, n, \mu} \, \mcV_{l, \mu}(\bfx) \, P_{n}^{(\beta , l)}(\rho) \Big) \, d\Omega \\
& \quad = \ 
\big( (n + 1)^{\gamma - s} \, (n + l + 1)^{\epsilon - t} \, p(\bfx) \, , \, (n + 1)^{s} \, (n + l + 1)^{t} q(\bfx) \big)_{L^{2}_{\beta}}  \, .
\end{align*}

In place of Lemma \ref{lma_nmeq1} we have the following.
\begin{corollary} \label{cor_nmeq1}
Let $\wtilde{u} \in \mcH^{1 + s}_{\alpot}(\Omega)$ and 
$U \in \big( H^{s}_{\alpot - 1}(\Omega) \big)^{2}$ be given by
$\Grad \wtilde{u} \, = \, \omega^{\alpot - 1} \, U$. Then
\be
     \frac{9  \alpha \, + \, 2}{8}  \, \, \| \wtilde{u} \|^{2}_{\mcH^{1 + s}_{\alpot}(\Omega)} \ \le \ 
     \| U \|^{2}_{\big( H^{s}_{\alpot - 1}(\Omega) \big)^{2}}
     \ \le \  4 \,  \| \wtilde{u} \|^{2}_{\mcH^{1 + s}_{\alpot}(\Omega)} \, .
\label{HHnmeq1} 
\ee
\end{corollary}
\textbf{Proof}: The proof is almost verbatim that of Lemma \ref{lma_nmeq1}, with the additional weights $(n + 1)^{s} \, (n + l + 1)^{s}$
cancelling out in equations \eqref{rato1} - \eqref{rato4}. \\
\mbox{ } \hfill \qed

In place of Lemma \ref{lmau2v} we have the following.
\begin{corollary} \label{cor_au2v}
Let $\wtilde{u}(\bfx) \ = \ \omega^{\alpot} \, \sum_{l, n, \mu} a_{l, n, \mu} \, \mcV_{l, \mu}(\bfx) \, P_{n}^{(\alpot , l)}(\rho) \, \in \, 
\mcH_{\alpot}^{1 + s}(\Omega)$. Then
\begin{align}
& \wtilde{v}(\bfx) \ = \ 2^{2 - \alpha} \, \omega^{\alpot} \, \sum_{l, n, \mu} (n + 1) \, \frac{\Gamma(n + 1)}{\Gamma(n + 1 + \alpot)} \, 
 \frac{\Gamma(n + l + 1)}{\Gamma(n + l + \alpot)} \, 
 a_{l, n, \mu} \, \mcV_{l, \mu}(\bfx) \, P_{n}^{(\alpot , l)}(\rho) \, \in \, \mcH_{\alpot}^{\alpha - 1 + s}(\Omega) \, ,  \label{defvtda2} \\
 & \hspace{1.5in} \ (n + 1)^{s} \, (n + l + 1)^{s} \, \wtilde{v}(\bfx) \, \in \, \mcH_{\alpot}^{\alpha - 1 - s}(\Omega) \, ,  \nonumber \\
 & \hspace{1.0in} \ \mbox{with } \ \| \wtilde{u} \|_{\mcH_{\alpot}^{1 + s}(\Omega)} \ 
 \sim  \ \|  (n + 1)^{s} \, (n + l + 1)^{s} \, \wtilde{v} \|_{\mcH_{\alpot}^{\alpha - 1 - s}(\Omega)}  \, . \nonumber 
\end{align}
\end{corollary}
\textbf{Proof}: The proof is almost verbatim that of Lemma \ref{lmau2v}. \\
\mbox{ } \hfill \qed

In place of Lemma \ref{lmaBdDif} we have the following.
\begin{corollary} \label{cor_BdDif}
Let $\wtilde{u}(\bfx)$ and $\wtilde{v}(\bfx)$ be defined as in Lemma \ref{lmau2v}, and $U$ and $V$ given by
$\Grad \wtilde{u} \ = \ \omega^{\alpot - 1} \, U$ and $( - \Delta )^{\frac{2 - \alpha}{2}} \Grad \wtilde{v} \ = \ V$.
Then for $W \ = \ U \, - \, V$ we have
\begin{align}
\| (n + 1)^{\frac{s}{2}} \, (n + l + 1)^{\frac{s}{2}} \, W \|_{L^{2}_{\alpot - 1}(\Omega)} \ 
&\le \ \frac{(2 - \alpha)}{\sqrt{\alpha \, (2 + \alpha)}} \, 
\| (n + 1)^{\frac{s}{2}} \, (n + l + 1)^{\frac{s}{2}} \, U \|_{L^{2}_{\alpot - 1}(\Omega)}    \nonumber  \\
&= \  \frac{(2 - \alpha)}{\sqrt{\alpha \, (2 + \alpha)}} \, \| U \|_{H^{s}_{\alpot - 1}(\Omega)}  \, . 
\label{WWbd1}
\end{align}
\end{corollary}
\textbf{Proof}: The proof follows analogously to that for Lemma \ref{lmaBdDif}. The weights multiplying $W$ and $U$
cancel out in the equations \eqref{dfx1} - \eqref{dfx2}. \\
\mbox{ } \hfill \qed

The higher regularity result for $\wtilde{u}$ is given in the following corollary.
\begin{corollary}  \label{corexun2}
Assume that $1 < \alpha < 2$, $s > -1$ and $K \in \big( W_{w}^{\ceil{ |s| }, \infty}(\Omega) \big)^{2 \times 2}$ 
a symmetric positive definite matrix
with $\lambda_{m}$ and $\lambda_{M}$ in \eqref{Kspd} satisfying 
$\lambda_{M} \, < \, \frac{\sqrt{\alpha \, (2 + \alpha)}}{(2 - \alpha)} \, \lambda_{m}$.
Then, given $f \in H^{- (\alpha - 1) \, + \, s}_{\alpot}(\Omega)$ there exists a unique 
$\wtilde{u} \in \mcH^{1 + s}_{\alpot}(\Omega)$ satisfying
\eqref{sseq1}, with $\| \wtilde{u} \|_{\mcH^{1 + s}_{\alpot}(\Omega)} \ \le \ \frac{1}{C_{2}} \, \| f \|_{H^{-(\alpha - 1) + s}_{\alpot}(\Omega)}$,
where $C_{2}$ is the $inf-sup$ constant (see \eqref{condd2}).
\end{corollary}
\textbf{Proof}: 
To establish the continuity of $B(\cdot , \cdot)$ in \eqref{sseq1}, in order to avoid the issue of interpreting the product of
a distribution and a function, we consider two cases. For $s \ge 0$ we associate $K(\bfx)$ with $U$ and for the case
$-1 < s < 0$ we associate $K(\bfx)$ with the $( - \Delta )^{\frac{\alpha - 2}{2}} \Grad \, \wtilde{v}$ terms.

For $s \ge 0$, using \eqref{wsp3}, \eqref{HHnmeq1} and Corollary \ref{cormapdivG},   we have
\begin{align*}
K(\bfx) \, U \ \in H^{s}_{\alpot - 1}(\Omega) , & \mbox{  with }
\| K(\bfx) \, U \|_{H^{s}_{\alpot - 1}} \ \le \ \| K \|_{W_{w}^{\ceil{s}, \infty}} \, \| U \|_{H^{s}_{\alpot - 1}} \ 
\sim \  \| K \|_{W_{w}^{\ceil{s}, \infty}} \, \| \wtilde{u} \|_{\mcH^{1 + s}_{\alpot}} ,  \\
& \mbox{ and  }
\| ( - \Delta )^{\frac{\alpha - 2}{2}} \Grad \, \wtilde{v} \|_{H^{-s}_{\alpot - 1}} \ \sim \ \, \| \wtilde{v} \|_{\mcH^{\alpha - 1 - s}_{\alpot}} \, . 
\end{align*}
\[
\mbox{Thus,  } \ 
| B(\wtilde{u} \, , \, \wtilde{v}) | \  \le \ \| K(\bfx) \, U \|_{H^{s}_{\alpot - 1}}  \, 
\| ( - \Delta )^{\frac{\alpha - 2}{2}} \Grad \, \wtilde{v} \|_{H^{-s}_{\alpot - 1}} 
\ \le \ \| \wtilde{u} \|_{\mcH^{1 + s}_{\alpot}} \, \| \wtilde{v} \|_{\mcH^{\alpha - 1 - s}_{\alpot}} \, .
\]
A similar argument (with $K(\bfx)$ associated with $( - \Delta )^{\frac{\alpha - 2}{2}} \Grad \, \wtilde{v}$) establishes the continuity of
$B(\cdot , \cdot)$ for $-1 < s < 0$.

To establish the $inf - sup$ condition, for $\wtilde{u} \in \mcH^{1 + s}_{\alpot}(\Omega)$ and
$\wtilde{v} \in \mcH^{\alpha - 1 + s}_{\alpot}(\Omega)$ given by \eqref{defvtda2}, consider
\begin{align*}
 & \sup_{ \wtilde{z} \in \mcH_{\alpot}^{\alpha - 1 - s}} 
   \frac{ | B(\wtilde{u} \, , \, \wtilde{z})  |}{ \| \wtilde{z} \|_{ \mcH_{\alpot}^{\alpha - 1 - s}}}  \ \ge \
    \frac{ | B(\wtilde{u} \ , \ (n + 1)^{s} \, (n + l + 1)^{s} \, \wtilde{v})  |}%
    { \| (n + 1)^{s} \, (n + l + 1)^{s} \, \wtilde{v} \|_{ \mcH_{\alpot}^{\alpha - 1 - s}}} \  \\
 & \hspace*{1.4in}
    \sim \ \frac{ | B(\wtilde{u} \ , \ (n + 1)^{s} \, (n + l + 1)^{s} \, \wtilde{v})  |}{ \| \wtilde{u} \|_{ \mcH_{\alpot}^{1 + s}}}  
    \ \mbox{ (using Corollary \ref{cor_au2v})}  \\
&= \ \frac{1}{ \| \wtilde{u} \|_{ \mcH_{\alpot}^{1 + s}} } \, \int_{\Omega} K \, \omega^{\alpot - 1} \, U \cdot 
(n + 1)^{s} \, (n + l + 1)^{s} \, V \, d\Omega   \\
&= \ \frac{1}{ \| \wtilde{u} \|_{ \mcH_{\alpot}^{1 + s}} } \, 
\Big( \int_{\Omega} K \, \omega^{\alpot - 1} \, U \cdot (n + 1)^{s} \, (n + l + 1)^{s} \, U \, d\Omega  
\ - \ \int_{\Omega} K \, \omega^{\alpot - 1} \, U \cdot (n + 1)^{s} \, (n + l + 1)^{s} \, W \, d\Omega  \Big)  \\
&\ge \ \frac{1}{ \| \wtilde{u} \|_{ \mcH_{\alpot}^{1 + s}} } \, 
\Big( \lambda_{m} \, \| (n + 1)^{s / 2} \, (n + l + 1)^{s / 2} \,U \|_{L^{2}_{\alpot - 1}}^{2}   \\
& \hspace*{1.5in}
\ - \ \lambda_{M} \, \| (n + 1)^{s / 2} \, (n + l + 1)^{s / 2} \, U \|_{L^{2}_{\alpot - 1}} \, 
\| (n + 1)^{s / 2} \, (n + l + 1)^{s / 2} \, W \|_{L^{2}_{\alpot - 1}} \Big) \\
&\ge \ \frac{1}{ \| \wtilde{u} \|_{ \mcH_{\alpot}^{1 + s}} } \, \big( \lambda_{m} \ - \ 
\frac{(2 - \alpha)}{\sqrt{\alpha \, (2 + \alpha)}}  \, \lambda_{M} \big)
\, \| (n + 1)^{s / 2} \, (n + l + 1)^{s / 2} \, U \|_{L^{2}_{\alpot - 1}}^{2}  \ \mbox{ (using Corollary \ref{cor_BdDif})} \\
&\ge \ \frac{9 \alpha \, + \, 2}{8} \, \big( \lambda_{m} \ - \ 
\frac{(2 - \alpha)}{\sqrt{\alpha \, (2 + \alpha)}}  \, \lambda_{M} \big) \, \| \wtilde{u} \|_{ \mcH_{\alpot}^{1 + s}}
 \ \mbox{ (using \eqref{HHnmeq1})}  \\
&= \ C_{2} \,  \| \wtilde{u} \|_{ \mcH_{\alpot}^{1 + s}} \, , \ \mbox{ provided } 
\lambda_{M} \, < \, \frac{\sqrt{\alpha \, (2 + \alpha)}}{(2 - \alpha)}  \, \lambda_{m} \, .   
\end{align*}

The condition
\[
  \sup_{\wtilde{u} \in \mcH_{\alpot}^{1 + s}(\Omega)} 
    | B(\wtilde{u} \, , \, \wtilde{z})  |   \ > \ 
     0 \, , \ \ \forall  \ \ 0 \neq \wtilde{z} \in \mcH_{\alpot}^{\alpha - 1 - s}(\Omega) \, ,
\] 
follows using a similar argument as used for the $inf - sup$ condition.

The norm bound for the solution $\wtilde{u}$ is obtained as in the proof of Theorem \ref{thmexun}. \\
\mbox{ } \hfill \qed

 \setcounter{equation}{0}
\setcounter{figure}{0}
\setcounter{table}{0}
\setcounter{theorem}{0}
\setcounter{lemma}{0}
\setcounter{corollary}{0}
\setcounter{definition}{0}
\section{Concluding Remarks}
\label{sec_conc}
In this paper we have analyzed the existence, uniqueness and regularity of the solution to the generalized, 
variable diffusivity,
fractional Laplace equation on the unit disk in $\real^{2}$. Our results show that for the symmetric, positive definite, diffusivity
matrix, $K(\bfx)$, satisfying $\lambda_{m} \, \bfv^{T} \, \bfv \ \le \, \bfv^{T} \, K(\bfx) \, \bfv \ \le \ \lambda_{M} \, \bfv^{T} \, \bfv$, 
for all $\bfv \in \real^{2}$, $\bfx \in \Omega$, with 
$\lambda_{M} \, < \, \frac{\sqrt{\alpha \, (2 + \alpha)}}{(2 - \alpha)} \, \lambda_{m}$, the 
problem has a unique solution. The regularity of the solution is given in an appropriately weighted Sobolev space in 
terms of the regularity of the right hand side function and $K(\bfx)$.

It is interesting to note that as $\alpha \rightarrow 2$, in which case the problem becomes the usual variable coefficient
Laplace equation, the condition on $\lambda_{M}$ disappears. This is consistent with the analysis for the usual variable coefficient
Laplace equation.

In \cite{zhe232} existence and uniqueness of the generalized, 
fractional Laplace equation on the unit disk in $\real^{2}$ was shown,
for $K(\bfx)$ a constant, symmetric, positive definite matrix. No restriction on $\lambda_{M}$ was needed to establish the
existence and uniqueness of solution.

It remains an open question for the general case of $K(\bfx)$ if the sufficient condition 
$\lambda_{M} \, < \, \frac{\sqrt{\alpha \, (2 + \alpha)}}{(2 - \alpha)} \, \lambda_{m}$ is also necessary. \\

\textbf{Acknowledgement}: The author gratefully acknowledges helpful discussions with 
Professors Hong Wang (University of South Carolina, U.S.A.) and Xiangcheng Zhang (Shandong University, China).


\appendix
 
 \setcounter{equation}{0}
\setcounter{figure}{0}
\setcounter{table}{0}
\setcounter{theorem}{0}
\setcounter{lemma}{0}
\setcounter{corollary}{0}
\setcounter{definition}{0}
\section{Proofs of the mapping properties}
\label{sec_mapPrf}
In this section we present the proofs of the mapping properties described by Corollary \ref{cormapD2} and Corollary \ref{cormapG}.
We begin by giving some properties of the norms of the basis functions.

Let $C_{l, \mu} \ = \ \left\{ \begin{array}{l}
                                  \pi/2 , \ \mbox{if } l = 0, \ \mu = 1,  \\
                                  \pi , \ \mbox{for } l \ge 1, 
                                  \end{array}  \right. $
\begin{lemma} \cite[Lemma B.4]{erv241}   \label{lmaeq3}
For $l, \, n, \, k \ge 0$, $(l, \, \mu) \neq (0, -1)$,
\[ 
\frac{\| \mcV_{l, \mu}(\bfx) \, P_{n}^{(\gamma + k \, , \, l)}(\rho) \|^{2}_{L^{2}_{\gamma + k}}}%
{  \| \mcV_{l, \mu}(\bfx) \, P_{n}^{(\gamma , l)}(\rho) \|^{2}_{L^{2}_{\gamma}} }
  \ = \ \frac{2n + \gamma + l + 1}{2n + \gamma + l + 1 + k} \, 
 \frac{ \Pi_{s = 1}^{k} \, (n + \gamma + s)}{ \Pi_{s = 1}^{k} \, (n + \gamma + l + s)} \,  .
\]
\end{lemma}

\begin{lemma}  \cite[Lemma B.5]{erv241}   \label{lmaeq4}  
For $l, \, n, \, j, \, m \geq 0$, $(l, \, \mu) \neq (0, -1)$,
\begin{align*}
& \frac{ \| \mcV_{l+j \, , \, \mu}(\bfx) \, P_{n+m}^{(\gamma , \, l + j )}(\rho) \|^{2}_{L^{2}_{\gamma}} }%
{ \| \mcV_{l, \mu}(\bfx) \, P_{n}^{(\gamma , l)}(\rho) \|^{2}_{L^{2}_{\gamma}} }   \\   
& \quad \quad \quad \ = \ \frac{C_{l+j , \mu}}{C_{l , \mu}} \, \frac{2n + \gamma + l + 1}{2n + 2m + \gamma + l + j + 1} \, 
 \frac{ \Pi_{s = 1}^{m} \, (n + \gamma + s)}{ \Pi_{s = 1}^{m} \, (n +  s)} \, 
 \frac{ \Pi_{s = 1}^{m+j} \, (n + l + s)}{ \Pi_{s = 1}^{m+j} \, (n + \gamma + l + s)} \, 
 \, .
\end{align*}
\end{lemma}

\begin{lemma}  \cite[Lemma B.6]{erv241}   \label{lmaeq5}
For $l, \, n \ge 0$ and $m \geq j \geq 0$ with $(l - j , \, \mu) \neq (0, -1)$,
\begin{align*} 
& \frac{ \| \mcV_{l-j \, , \, \mu}(\bfx) \, P_{n+m}^{(\gamma , \, l - j )}(\rho) \|^{2}_{L^{2}_{\gamma}} }%
{ \| \mcV_{l, \mu}(\bfx) \, P_{n}^{(\gamma , l)}(\rho) \|^{2}_{L^{2}_{\gamma}} }   \\    
& \quad \quad \quad \ = \ \frac{C_{l-j , \mu}}{C_{l , \mu}} \, \frac{2n + \gamma + l + 1}{2n + 2m + \gamma + l - j + 1} \, 
 \frac{ \Pi_{s = 1}^{m} \, (n + \gamma + s)}{ \Pi_{s = 1}^{m} \, (n +  s)} \, 
 \frac{ \Pi_{s = 1}^{m-j} \, (n + l + s)}{ \Pi_{s = 1}^{m-j} \, (n + \gamma + l + s)} \, 
 \, .
\end{align*}
\end{lemma}

\textbf{Proof of Corollary \ref{cormapD2}} \\
Suppose 
\[
v(\bfx) \ = \ \sum_{l, n, \mu} a_{l, n, \mu} \, \mcV_{l, \mu}(\bfx) P_{n}^{(\gamma , l)}(\rho) \ \in H^{s}_{\gamma}(\Omega) \, ,
\]
i.e.,
\[
  \sum_{l, n, \mu} (n + 1)^{s} (n + l + 1)^{s} \, 
  a_{l, n, \mu}^{2} \, \| \mcV_{l, \mu}(\bfx) P_{n}^{(\gamma , l)}(\rho) \|_{L^{2}_{\gamma}}^{2} \ < \ \infty \, .
\]
From \eqref{dersw1},  keeping in mind that $a_{l, n, \mu} = 0$ for $l, n < 0$, $a_{0, n, -1} = 0$,  $\mcV_{0, -1}(\bfx) = 0$,
(after reindexing)
\[
\frac{\partial v(\bfx)}{\partial x} \ = \ 
 \sum_{l, n, \mu} \Big( (n + l + 1) \, a_{l+1, n, \mu} \ + \ (n + \gamma + l + 1) \, a_{l-1, n+1, \mu} \Big) \, 
 \mcV_{l, \mu}(\bfx) P_{n}^{(\gamma + 1 \,  , l)}(\rho) \, .
 \]
 Hence, using orthogonality,
 \begin{align*}
 \left\| \frac{\partial v(\bfx)}{\partial x}  \right\|^{2}_{H^{s-1}_{\gamma + 1}}
&\le \  \sum_{n} (n + 1)^{s-1} (n + 1)^{s-1} \, (n + 1)^{2} \, a_{1, n, 1}^{2} 
   \| \mcV_{0, 1}(\bfx) P_{n}^{(\gamma + 1 \,  , 0)}(\rho) \|_{L^{2}_{\gamma + 1}}^{2}    \\
&+ \ 2 \,  \sum_{l \ge 1, n, \mu} (n + 1)^{s-1} (n + l + 1)^{s-1} \, (n + l + 1)^{2} \, a_{l+1, n, \mu}^{2} 
   \| \mcV_{l, \mu}(\bfx) P_{n}^{(\gamma + 1 \,  , l)}(\rho) \|_{L^{2}_{\gamma + 1}}^{2}    \\
&+ \ 2 \,  \sum_{l \ge 1, n, \mu} (n + 1)^{s-1} (n + l + 1)^{s-1} \, (n + \gamma + l + 1)^{2} \, a_{l-1, n+1, \mu}^{2} 
   \| \mcV_{l, \mu}(\bfx) P_{n}^{(\gamma + 1 \,  , l)}(\rho) \|_{L^{2}_{\gamma + 1}}^{2}  \, . 
\end{align*}
 
Using Lemma \ref{lmaeq3} (with $k = 1$), and Lemma \ref{lmaeq4} (with $j = 1$, $m = 0$),
\begin{align}
& \| \mcV_{l, \mu}(\bfx) P_{n}^{(\gamma + 1 \,  , l)}(\rho) \|_{L^{2}_{\gamma + 1}}^{2}
\ = \ \frac{2n + \gamma + l + 1}{2n + \gamma + l + 2} \, \frac{(n + \gamma + 1)}{(n + \gamma + l + 1)} \,
\| \mcV_{l, \mu}(\bfx) P_{n}^{(\gamma  , l)}(\rho) \|_{L^{2}_{\gamma}}^{2}   \label{kilo1}   \\
&= \ \frac{2n + \gamma + l + 1}{2n + \gamma + l + 2} \, \frac{(n + \gamma + 1)}{(n + \gamma + l + 1)} \,
    \frac{C_{l, \mu}}{C_{l+1, \mu}} \, \frac{2n + \gamma + l + 2}{2n + \gamma + l + 1} \, \frac{(n + \gamma + l + 1)}{(n + l + 1)} 
 \| \mcV_{l+1, \mu}(\bfx) P_{n}^{(\gamma  , l+1)}(\rho) \|_{L^{2}_{\gamma}}^{2}   \nonumber \\
&= \ \frac{C_{l, \mu}}{C_{l+1, \mu}} \, \frac{(n + \gamma + 1)}{(n + l + 1)} \,
  \| \mcV_{l+1, \mu}(\bfx) P_{n}^{(\gamma  , l+1)}(\rho) \|_{L^{2}_{\gamma}}^{2}  \, .  \label{kilo2}
\end{align}
 
From \eqref{kilo1} and Lemma \ref{lmaeq5} (with $j = 1$, $m = 1$),
\begin{align}
& \| \mcV_{l, \mu}(\bfx) P_{n}^{(\gamma + 1 \,  , l)}(\rho) \|_{L^{2}_{\gamma + 1}}^{2}   \nonumber \\
&= \ \frac{2n + \gamma + l + 1}{2n + \gamma + l + 2} \, \frac{(n + \gamma + 1)}{(n + \gamma + l + 1)} \,
    \frac{C_{l, \mu}}{C_{l-1, \mu}} \, \frac{2n + \gamma + l + 2}{2n + \gamma + l + 1} \, \frac{(n + 1)}{(n + \gamma + 1)} 
 \| \mcV_{l-1, \mu}(\bfx) P_{n+1}^{(\gamma  , l-1)}(\rho) \|_{L^{2}_{\gamma}}^{2}   \nonumber \\
&= \  \frac{C_{l, \mu}}{C_{l-1, \mu}} \, \frac{(n + 1)}{(n + \gamma + l + 1)} \,
  \| \mcV_{l-1, \mu}(\bfx) P_{n+1}^{(\gamma  , l-1)}(\rho) \|_{L^{2}_{\gamma}}^{2}  \, .  \label{kilo3}
\end{align}

With \eqref{kilo2} and \eqref{kilo3} we obtain
\begin{align*}
& \left\| \frac{\partial v(\bfx)}{\partial x}  \right\|^{2}_{H^{s-1}_{\gamma + 1}}
\ \lesssim \\
& \sum_{n} (n + 1)^{s-1} (n + 1)^{s-1} \, (n + 1) \, (n + \gamma + 1) a_{1, n, 1}^{2} 
   \| \mcV_{1, 1}(\bfx) P_{n}^{(\gamma + 1 \,  , 1)}(\rho) \|_{L^{2}_{\gamma}}^{2}    \\
&+  \ \sum_{l \ge 1, n, \mu} (n + 1)^{s-1} (n + l + 1)^{s+1} \, \frac{(n + \gamma + 1)}{(n + l + 1)} \, a_{l+1, n, \mu}^{2}  \, 
   \| \mcV_{l+1, \mu}(\bfx) P_{n}^{(\gamma  , l+1)}(\rho) \|_{L^{2}_{\gamma}}^{2}    \\
&+ \   \sum_{l \ge 1, n, \mu} (n + 1)^{s-1} (n + l + 1)^{s-1} \, (n + \gamma + l + 1)^{2} \, \frac{(n + 1)}{(n + \gamma + l + 1)} \, 
a_{l-1, n+1, \mu}^{2}  \,   \| \mcV_{l-1, \mu}(\bfx) P_{n+1}^{(\gamma  , l-1)}(\rho) \|_{L^{2}_{\gamma}}^{2} \\
%
&\lesssim \ \sum_{l, n, \mu} (n + 1)^{s} (n + l + 1)^{s}  \, a_{l, n, \mu}^{2} \, 
       \| \mcV_{l, \mu}(\bfx) P_{n}^{(\gamma  , l)}(\rho) \|_{L^{2}_{\gamma}}^{2}   \\
&= \ \| v \|_{H^{s}_{\gamma}}^{2} \ < \ \infty \, ,
\end{align*}
which establishes the bounded mapping property for $\frac{\partial}{\partial x}$. The mapping property for $\frac{\partial}{\partial y}$ 
is established in an analogous manner.  \\
\mbox{  } \hfill \qed

\textbf{Proof of Corollary \ref{cormapG}} \\
Let $\wtilde{v}(\bfx) \ =  \ \omega^{\gamma} v(\bfx)$, where
\[
v(\bfx) \ = \ \sum_{l, n, \mu} a_{l, n, \mu} \, \mcV_{l, \mu}(\bfx) P_{n}^{(\gamma , l)}(\rho) \ \in H^{s}_{\gamma}(\Omega) \, ,
\]
i.e.,
\[
  \sum_{l, n, \mu} (n + 1)^{s} (n + l + 1)^{s} \, 
  a_{l, n, \mu}^{2} \, \| \mcV_{l, \mu}(\bfx) P_{n}^{(\gamma , l)}(\rho) \|_{L^{2}_{\gamma}}^{2} \ < \ \infty \, .
\]

From \eqref{ders1}, keeping in mind that $a_{l, n, \mu} = 0$ for $l, n < 0$, $a_{0, n, -1} = 0$, $\mcV_{0, -1}(\bfx) = 0$,
\begin{align*}
\frac{\partial \wtilde{v}}{\partial x} &= \ 
- \omega^{\gamma - 1} \Big\{ \sum_{l, n, \mu} (n + \gamma) a_{l, n, \mu} \, \mcV_{l+1 \, , \mu}(\bfx) P_{n}^{(\gamma-1 \,  , \,  l+1)}(\rho)
   \\
& \quad \quad \quad   \quad \quad \quad
 \ + \ (n + 1) a_{l, n, \mu} \, \mcV_{l-1 \, , \mu}(\bfx) P_{n+1}^{(\gamma-1 \,  , \,  l-1)}(\rho) \Big\}  \\
&:= - \omega^{\gamma - 1} \, w_{1}(\bfx) \,   \\
\mbox{where, } \ w_{1}(\bfx) &= \ \sum_{l, n, \mu} \left( (n + \gamma) \, a_{l-1, n, \mu} \ + \ n \, a_{l+1, n-1, \mu} \right) 
\mcV_{l, \mu}(\bfx) P_{n}^{(\gamma-1 , l)}(\rho)  \, ,
\end{align*}
and from \eqref{ders2},
\begin{align*}
\frac{\partial \wtilde{v}}{\partial y} &= \ 
- \omega^{\gamma - 1} \Big\{ \sum_{l, n, \mu} (\pm) (n + \gamma) a_{l, n, \mu} \, \mcV_{l+1 \, , \mu^{*}}(\bfx) 
P_{n}^{(\gamma-1 \,  , \,  l+1)}(\rho)
   \\
& \quad \quad \quad   \quad \quad \quad
 \ + \  (\mp) (n + 1) a_{l, n, \mu} \, \mcV_{l-1 \, , \mu^{*}}(\bfx) P_{n+1}^{(\gamma-1 \,  , \,  l-1)}(\rho) \Big\}  \\
&:= - \omega^{\gamma - 1} \, w_{2}(\bfx) \,   \\
\mbox{where, } \ w_{2}(\bfx) &= \ \sum_{l, n, \mu} \left( (n + \gamma) \, a_{l-1, n, \mu} \ - \ n \, a_{l+1, n-1, \mu} \right) 
(\pm) \mcV_{l, \mu^{*}}(\bfx) P_{n}^{(\gamma-1 , l)}(\rho)  \, .
\end{align*}

Next,
\begin{align}
& \| w_{1} \|^{2}_{H^{s-1}_{\gamma - 1}} \ + \  \| w_{2} \|^{2}_{H^{s-1}_{\gamma - 1}}  \nonumber \\
& \ = \ 
 \sum_{n \ge 1} (n + 1)^{s-1} \, (n + 1)^{s-1} \, n^{2} \, \left(  a_{1, n-1, 1}^{2} \ + \ a_{1, n-1, -1}^{2} \right) \, 
 \|  \mcV_{0, 1}(\bfx) P_{n}^{(\gamma-1 , 0)}(\rho) \|^{2}_{L^{2}_{\gamma - 1}}  \nonumber \\
&\quad + \  \sum_{l \ge 1, n, \mu} (n + 1)^{s-1} \, (n + l + 1)^{s-1} \, \left( (n + \gamma) \, 
a_{l-1, n, \mu} \ + \ n \, a_{l+1, n-1, \mu} \right)^{2} \, 
 \|  \mcV_{l, 1}(\bfx) P_{n}^{(\gamma-1 , l)}(\rho) \|^{2}_{L^{2}_{\gamma - 1}}  \nonumber \\
& \quad + \ 
 \sum_{l \ge 1, n, \mu} (n + 1)^{s-1} \, (n + l + 1)^{s-1} \, \left( (n + \gamma) \, a_{l-1, n, \mu} \ - \ n \, a_{l+1, n-1, \mu} \right)^{2} \, 
 \|  \mcV_{l, 1}(\bfx) P_{n}^{(\gamma-1 , l)}(\rho) \|^{2}_{L^{2}_{\gamma - 1}}  \nonumber \\
& \ = \ 
 \sum_{n \ge 1} (n + 1)^{s-1} \, (n + l + 1)^{s-1} \, n^{2} \, \left(  a_{1, n-1, 1}^{2} \ + \ a_{1, n-1, -1}^{2} \right) \, 
 \|  \mcV_{0, 1}(\bfx) P_{n}^{(\gamma-1 , 0)}(\rho) \|^{2}_{L^{2}_{\gamma - 1}}  \nonumber \\
&   \quad + \ 
2 \, \sum_{l \ge 1, n, \mu} (n + 1)^{s-1} \, (n + l + 1)^{s-1} \, \left( (n + \gamma)^{2} \, a_{l-1, n, \mu}^{2} \ 
 + \ n^{2} \, a_{l+1, n-1, \mu}^{2} \right)^{2} \, 
 \|  \mcV_{l, 1}(\bfx) P_{n}^{(\gamma-1 , l)}(\rho) \|^{2}_{L^{2}_{\gamma - 1}}  \, .  \label{poik1}
\end{align} 

From Lemma \ref{lmaeq3} (with $\gamma \rightarrow \gamma-1$, $k = 1$),
\be
 \|  \mcV_{l, \mu}(\bfx) P_{n}^{(\gamma-1 , l)}(\rho) \|^{2}_{L^{2}_{\gamma - 1}} \ = \
 \frac{2n + \gamma + l + 1}{2n + \gamma + l} \, \frac{n + \gamma + l}{n + \gamma} \,  
  \|  \mcV_{l, \mu}(\bfx) P_{n}^{(\gamma , l)}(\rho) \|^{2}_{L^{2}_{\gamma}} \, .
  \label{poik2}
\ee

Using Lemma \ref{lmaeq4} (with $l \rightarrow l-1$, $j = 1$, $m = 0$)  and \eqref{poik2},
\begin{align}
 \|  \mcV_{l, \mu}(\bfx) P_{n}^{(\gamma-1 , l)}(\rho) \|^{2}_{L^{2}_{\gamma - 1}} &= \
 \frac{2n + \gamma + l + 1}{2n + \gamma + l} \, \frac{n + \gamma + l}{n + \gamma} \,   \nonumber  \\
 & \quad \quad \cdot \ \frac{C_{l, \mu}}{C_{l-1, \mu}} \, \frac{2n + \gamma + l}{2n + \gamma + l + 1} \, 
 \frac{n + l}{n + \gamma + l} \,  \|  \mcV_{l-1, \mu}(\bfx) P_{n}^{(\gamma , l-1)}(\rho) \|^{2}_{L^{2}_{\gamma}} \nonumber \\
 &= \ \   \frac{C_{l, \mu}}{C_{l-1, \mu}} \, 
   \frac{n + l}{n + \gamma} \,  \|  \mcV_{l-1, \mu}(\bfx) P_{n}^{(\gamma , l-1)}(\rho) \|^{2}_{L^{2}_{\gamma}} \, .
  \label{poik3}
\end{align}

Using Lemma \ref{lmaeq5} (with $l \rightarrow l+1$, $j \rightarrow 1$, $n \rightarrow n-1$, $m \rightarrow 1$)
and \eqref{poik2},
\begin{align}
 \|  \mcV_{l, \mu}(\bfx) P_{n}^{(\gamma-1 , l)}(\rho) \|^{2}_{L^{2}_{\gamma - 1}} &= \
 \frac{2n + \gamma + l + 1}{2n + \gamma + l} \, \frac{n + \gamma + l}{n + \gamma} \,   \nonumber  \\
 & \quad \quad \cdot \ \frac{C_{l, \mu}}{C_{l+1, \mu}} \, \frac{2n + \gamma + l}{2n + \gamma + l + 1} \, 
 \frac{n + \gamma}{n} \,  \|  \mcV_{l+1, \mu}(\bfx) P_{n-1}^{(\gamma , \, l+1)}(\rho) \|^{2}_{L^{2}_{\gamma}} \nonumber \\
&= \ \   \frac{C_{l, \mu}}{C_{l+1, \mu}} \, 
   \frac{n + \gamma + l}{n} \,  \|  \mcV_{l+1, \mu}(\bfx) P_{n-1}^{(\gamma , l+1)}(\rho) \|^{2}_{L^{2}_{\gamma}} \, .
  \label{poik4}
\end{align}

Combining \eqref{poik1} with \eqref{poik3} and \eqref{poik4}
\begin{align*}
&\| w_{1} \|^{2}_{H^{s-1}_{\gamma - 1}} \ + \ \| w_{2} \|^{2}_{H^{s-1}_{\gamma - 1}}  \\  
& \ = \ 
 \sum_{n \ge 1} (n + 1)^{s-1} \, (n + 1)^{s-1} \, n^{2} \, \left(  a_{1, n-1, 1}^{2} \ + \ a_{1, n-1, -1}^{2} \right) \, 
  \frac{C_{0, 1}}{C_{1, 1}} \, 
   \frac{n + \gamma}{n} \,  \|  \mcV_{1, 1}(\bfx) P_{n-1}^{(\gamma , 1)}(\rho) \|^{2}_{L^{2}_{\gamma}}   \\  
& \ \  + \ 
2 \, \sum_{l \ge 1, n, \mu} (n + 1)^{s-1} \, (n + l + 1)^{s-1} \, (n + \gamma)^{2} \, a_{l-1, n, \mu}^{2} \, 
  \frac{C_{l, \mu}}{C_{l-1, \mu}} \, 
   \frac{n + l}{n + \gamma} \,  \|  \mcV_{l-1, 1}(\bfx) P_{n}^{(\gamma , l-1)}(\rho) \|^{2}_{L^{2}_{\gamma}}  \\  
& \ \  + \   
 2 \,  \sum_{l \ge 1, n, \mu} (n + 1)^{s-1} \, (n + l + 1)^{s-1} \,
  n^{2} \, a_{l+1, n-1, \mu}^{2} \, 
 \frac{C_{l, \mu}}{C_{l+1, \mu}} \, 
   \frac{n + \gamma + l}{n} \,  \|  \mcV_{l+1, 1}(\bfx) P_{n-1}^{(\gamma , l+1)}(\rho) \|^{2}_{L^{2}_{\gamma}} \\  
&\sim \
  \sum_{l, n, \mu} (n + 1)^{s} \, (n + l + 1)^{s} \,   a_{l, n, \mu}^{2} \, 
    \|  \mcV_{l, \mu}(\bfx) P_{n}^{(\gamma , l)}(\rho) \|^{2}_{L^{2}_{\gamma}}  \, ,   \\  
&= \ \| v \|^{2}_{H^{s}_{\gamma}} \, .     
\end{align*}
Thus, $\Grad \wtilde{v} \, = \, \Grad \omega^{\gamma} v \, = \, \omega^{\gamma - 1} \left[ \begin{array}{c} w_{1} \\ w_{2} \end{array} \right]
\, := \, \omega^{\gamma - 1} \bfw \, \in \, \omega^{\gamma - 1} \otimes \left( H^{s-1}_{\gamma - 1}(\Omega) \right)^{2} $,
with $\| v \|_{H_{\gamma}^{s}} \sim \| \bfw \|_{H_{\gamma - 1}^{s-1}}$.

The norm equivalence of $v$ and $\bfw$ establishes the mapping is bounded and one-to one, and hence has a bounded inverse
on the image space of the mapping. \\
 \mbox{  } \hfill \qed



\end{document}